\documentclass [a4paper,twoside,12pt]{article}

\usepackage[utf8]{inputenc}
\usepackage{amsmath,amsfonts,amssymb}
\usepackage{vmargin,graphicx,theorem}
\usepackage[english]{babel}
\usepackage{enumerate}
\usepackage{color}

\RequirePackage[colorlinks,linkcolor=blue,citecolor=blue,urlcolor=blue]{hyperref}
\usepackage{amssymb, amsmath, amsfonts, latexsym}
\usepackage{enumerate,graphicx}
\usepackage{pst-fill,pst-grad,pst-plot,pst-eucl,pstricks-add,pst-node}

\setpapersize[portrait]{A4}
\setmarginsrb{1.5cm}{1cm}{1.5cm}{2.5cm}{1.5cm}{0cm}{0.5cm}{2cm}

\newcommand{\disp}{\displaystyle}

\newtheorem{theo}{Theorem}
\newtheorem{prop}[theo]{Proposition}
\newtheorem{lemma}[theo]{Lemma}
\newtheorem{cor}[theo]{Corollary}
\newtheorem{rem}[theo]{Remark}

\newtheorem{defi}[theo]{Definition}

\newcommand{\proofend}{~$\rhd$}
\newcommand{\proofbegin}{~$\lhd$}

\newenvironment{eproof}
               {\noindent {\emph{\textbf{Proof}}}\\\proofbegin~}
               {\proofend\\}

\newcommand{\ABS}[1]{\ensuremath{{\left| #1 \right|}}} 
\newcommand{\PAR}[1]{\ensuremath{{\left(#1\right)}}} 
\newcommand{\SBRA}[1]{\ensuremath{{\left[#1\right]}}} 
\newcommand{\BRA}[1]{\ensuremath{{\left\{#1\right\}}}} 

\renewcommand{\phi}{\varphi}

\renewcommand{\geq}{\geqslant}

\newcommand{\N}{\ensuremath{\mathbb{N}}} 
\newcommand{\R}{\ensuremath{\mathbb{R}}} 

\newcommand{\ent}{{\rm Ent}} 

\def\tr{\mathop{\rm tr}\nolimits} 

\newcommand{\upchi}{\raise1pt\hbox{$\chi$}}

\newcommand{\be}{\begin{equation}}
\newcommand{\ee}{\end{equation}}
\def\benu{\begin{enumerate}}
\def\eenu{\end{enumerate}}

\newcommand{\beq}{\begin{equation}}\newcommand{\eeq}{\end{equation}}
\begin{document}

\title{New sharp Gagliardo-Nirenberg-Sobolev inequalities and an improved Borell-Brascamp-Lieb inequality}
\author{François Bolley\thanks{Laboratoire de probabilit\'es et mod\`eles al\'eatoires,  Umr Cnrs 7599, Univ. Paris 6, francois.bolley@upmc.fr}, Dario Cordero-Erausquin\thanks{Institut de Math\'ematiques de Jussieu, Umr Cnrs 7586, Univ. Paris 6, France. dario.cordero@imj-prg.fr}, Yasuhiro Fujita\thanks{Toyama university, yfujita@sci.u-toyama.ac.jp},  Ivan Gentil\thanks{Univ Lyon, Univ; Claude Bernard Lyon 1, Umr Cnrs 5208, Institut Camille Jordan. gentil@math.univ-lyon1.fr}\\ and Arnaud Guillin\thanks{Laboratoire de Math\'ematiques Blaise Pascal, Umr Cnrs 6620, Univ.~Clermont Auvergne, France. guillin@math.univ-bpclermont.fr}}

\date{\today}

\maketitle

\abstract{We propose a new Borell-Brascamp-Lieb inequality which leads to novel sharp Euclidean inequalities such as Gagliardo-Nirenberg-Sobolev inequalities in $\R^n$ and in the half-space $\R^n_+$. This gives a new bridge between the geometric pont of view of the Brunn-Minkowski inequality and the functional point of view of the Sobolev type inequalities. In this way we  unify, simplify and results by S. Bobkov - M. Ledoux, M. del Pino - J. Dolbeault and B. Nazaret.  }

\bigskip

\medskip

\noindent
{\bf Mathematics Subject Classification (2000):}

\section{Introduction}

Sharp inequalities are interesting not only because they correspond to exact solution of variational problems (often related to problems in physics) but also because they encode in general deep geometric information on the underneath space. In the present paper, we shall be interested in a rather general new functional isoperimetric inequalities of Sobolev type, and their links with the Brunn-Minkowski inequality: 
\begin{equation}
\label{BM}
\mbox{vol}_n(A+B)^{1/n}\ge \mbox{vol}_n(A)^{1/n}+\mbox{vol}_n(B)^{1/n}
\end{equation}
for every non empty Borel bounded measurable sets $A$ in $\R^n,$
where $\mbox{vol}_n(\cdot)$ denotes Euclidean volume. If it is now classically known that sharp Sobolev inequalities (see e.g. \cite{bobkov-ledoux08}) may be derived through this Brunn-Minkowski inequality, we will see that via a new version of its functional counterpart, namely the Borell-Brascamp-Lieb inequality, we will be able to tackle both $\R^n$ case and half-space $\R^n_+$ case for sharp Sobolev and new Gagliardo-Nirenberg inequalities, in a rather simple and direct manner. In order to present this novel inequality, let us first introduce the general Sobolev inequalities in $\R^n$ which have inspired our line of thought.

\medskip

To simplify the notation, let $\|f\|_p=\|f\|_{L^{p}(\R^n)}$ denote the $L^p$-norm with respect to Lebesgue measure. The sharp classical Sobolev inequalities state that for $n \geq 2, p\in [1,n)$, $p^* = \frac{np}{n-p}$, and every smooth function $f$ on $\R^n$,
\be\label{classicalsob}
{\|f\|_p}\leq \frac{\|h_p\|_{p^*}}{\big(\int_{\R^n} |\nabla h_p |^p \big)^{1/p}}\PAR{\int_{\R^n} |\nabla f |^p}^{1/p}
\ee
with 
\be\label{classicalhp}
h_p(x):=(1+|x|^q)^{\frac{p-n}{p}}.
\ee
The optimal constants in the Sobolev inequalities have been first exhibited in~\cite{aubin76,talenti76}. Quite naturally, these inequalities admit a generalization when the Euclidean norm $|\cdot|$ on $\R^n$ is replaced by any norm or quasi-norm $\|\cdot\|$ on $\R^n$. Indeed, if we use a norm $\|\cdot\|$ to compute the size of the differential 
in~\eqref{classicalsob}, then the result remains true, 
\be
\label{classicalsob2}
{\|f\|_{p^*}}\leq \frac{\|h_p\|_{p^*}}{\big(\int_{\R^n} ||\nabla h_p ||_*^p \big)^{1/p}}\PAR{\int_{\R^n} ||\nabla f ||_*^p}^{1/p}
\ee
 where $\|y\|_\ast := \sup_{\|x\|\le 1}x\cdot y$. In this case, $h_p(x):=(1+||x||^q)^{\frac{p-n}{p}}$.

A natural generalization of this problem may then be the minimization, under integrability constraints on a function $g$, of more general quantities than $\int_{\R^n} ||\nabla g ||_*^p$, say of the form
$$
\int_{\R^n} F(\nabla g)\,  g^{\alpha} 
$$
where $F:\R^n \to \R\cup\{+\infty\}$ is a convex function. Note that we have to allow a term $g^\alpha$, $\alpha \in \R$ because it can no longer be put inside the gradient if $F$ is not homogeneous.

\medskip

A first answer in this direction, which is in fact an example of our main results, is the following in term of Sobolev type  inequality.

\begin{theo}[A first convex inequality]
\label{thm-IC2}
Let $n\geq 2$ and  $W:\R^n \to (0,+\infty)$ satisfying $\int W^{1-n}<+\infty$. Then for any smooth function $g$ such that  
$\int W^\ast \big( \nabla g \big)g^{-n}<+\infty$,  $\int g^{1-n}<+\infty$ and
$$
\int g^{-n} = \int W^{-n}=1
$$
one has
\begin{equation}
\label{eq-IC2}
\int W^\ast \big( \nabla g \big)g^{-n} \ge  \frac1{n-1} \int W^{1-n},
\end{equation}
with equality if $g$ is equal to $W$ and is convex.  \end{theo}

Here $W^\ast$ is the Legendre transform of the function $W$ (see below for details). 

\medskip

We shall see that the family of sharp Sobolev inequalities~\eqref{classicalsob2}, for $p\in [1, n)$ easily follows  from this theorem. Let us mention that the coefficients $n$ and $n-1$ in this theorem are not arbitrary at all: in some aspects, they are the ``good" ones to reach the Sobolev inequality, as we shall see. This may be compared to Corollary 2 of \cite{bobkov-ledoux08} which was derived via the Prekopa-Leindler inequality, leading to a more involved formulation and proof of the Sobolev inequalities.
\medskip

As mentioned above, our work is inspired by the Brunn-Minkowski-Borell theory. In turn, we are going to shed new light on this theory. It has been observed by S.~Bobkov and M. ~Ledoux in~\cite{bobkov-ledoux00,bobkov-ledoux08} that Sobolev inequalities can be reached through a functional version of the Brunn-Minkowski inequality, the so-called Borell-Brascamp-Lieb inequality, due to C.~Borell and H. J. Brascamp - E. H. Lieb (\cite{borell75,brascamp-lieb76}). However, one can not use the standard functional version of the inequality. Indeed, there is a subtle game with the dimension.  
The standard version states that, for $n\geq1$, given $s\in[0,1]$ (and $t=1-s$) and three positive functions 
$u,v,w:\R^n\to (0,\infty]$ such that $\int u = \int v=1$ and 
$$
\forall x,y\in \R^n,\qquad   w(sx+ ty) \ge (s\, u^{-1/n}(x) + t\, v^{-1/n}(y))^{-n},
$$
then 
$$\int w \ge  1.$$ Let us remark that we have stated here the ``strongest" version of Borell-Brascamp-Lieb inequality (say for the parameter $p=1/n$), see e.g. \cite[Th. 10.2]{gardner02}. 
By a simple change of functions, it turns out that the result can be stated as follows: 
let  three positive 
functions $g, W, H:\R^n\to (0,\infty]$ such that 
$$
\forall x,y\in \R^n,\qquad   H(sx+ ty) \le s\, g(x) + t \, W(y),
$$
and $\int W^{-n} = \int g^{-n}=1.$ Then
\begin{equation}
\label{eq-BBL}
\int H^{-n} \ge  1.
\end{equation}
In some sense, what is needed for the Sobolev inequality is to replace $n$ by the smaller ${n-1}$.  To do so, S. Bobkov and M. Ledoux  used a classical geometric strengthening of the Brunn-Minkowski inequality, for sets having an hyperplane section of same volume.

A natural question raised by S. Bobkov and M. Ledoux is whether the Sobolev inequality can be proved directly from a new kind of Borell-Brascamp-Lieb inequality. We will exhibit such a new functional inequality, that we believe is the correct one, in the sense that sharp (trace-) Sobolev inequalities (and actually the above Theorem~\ref{thm-IC2}) follow from it, and more generally  new (trace-) Gagliardo-Nirenberg inequalities; moreover it can be easily proved using a mass transportation argument. Its main form is the following:

\begin{theo}[An extended Borell-Brascamp-Lieb inequality]
\label{thm-BBL2}
Let $n\geq 2$. Let $g,W, H :\R^n \to (0,+\infty]$ be Borel  functions and $s\in [0,1]$ and $t=1-s$ be such that
$$
\forall x,y\in \R^n,\qquad  H(sx+ ty) \le s\, g(x) + t\, W(y), 
$$
and $\int W^{-n} = \int g^{-n}=1.$ Then
\begin{equation}
\label{eq-BBL2}
\int H^{1-n} \ge s\, \int g^{1-n} + t\, \int  W^{1-n} .
\end{equation}
\end{theo}

As we shall see, Theorem~\ref{thm-IC2} appears as a linearization of Theorem~\ref{thm-BBL2} for $t\rightarrow0$.

\medskip

There are more general families of Sobolev type inequalities that have attracted much attention these past years, namely the Gagliardo-Nirenberg inequalities in $\R^n$ of the form
$$
\|f\|_{\alpha}\leq C\,\|\nabla f \|_{p}^\theta  \, \|f\|_{\beta}^{1-\theta}.
$$ 
Here the coefficients $\alpha,\beta,p$ belong to some adequate range and $\theta\in[0,1]$ is fixed by a scaling condition. Sharp inequalities are known for a certain family of parameters since the pioneering work of M. del Pino and J. Dolbeault~\cite{delpino-dolbeault}: namely, for $p>1$, $\alpha=ap/(a-p)$ and $\beta=p(a-1)/(a-p)$ where $a>p$ is a free parameter.  

This family can be recovered by Theorem~\ref{thm-IC2}, or rather an extension of it (see Theorem~\ref{thm-case1}). In fact this extension turns out not only to be a natural way of recovering this family, but also allows to extend the family to parameters $a<p$ leading to the new Gagliardo-Nirenberg inequality  with negative powers
$$
\|f\|_{p\frac{a-1}{a-p}}\leq C\,\|\nabla f \|_{p}^\theta\,   \|f\|_{\frac{ap}{a-p}}^{1-\theta}.
$$ 
Here $p>a$ if $a\ge n+1$, or $p\in(a, \frac{n}{n+1-a})$ if $a\in[n,n+1)$, and $\theta$ is fixed by a scaling condition. Let us note that partial results for a narrower range of such $a<p$ have been proved by V.-H.~Nguyen~\cite{nguyen-sobolev}, by another approach. 

A crucial advantage of our approach is also its robustness: it can be applied to reach a new family of trace Gagliardo-Nirenberg-Sobolev inequalities which extend the 
trace  Sobolev inequality proved by B. Nazaret~\cite{nazaret}. Indeed, letting $\R_+^n=\{(u,x),\, u\geq0,\, x\in\R^{n-1}\}$ we obtain the sharp family of inequalities
$$
\|f\|_{L^\alpha(\partial\R_+^n) }\leq C\|\nabla f \|_{L^p(\R_+^n)}^\theta   \|f\|_{L^\beta(\R_+^n)}^{1-\theta}.
$$ 
Here $p>1$, $\alpha=p(a-1)/(a-p)$ and $\beta=p(a-1)/(a-p)$ where $a>p$ is a free parameter and again $\theta \in [0,1]$ is fixed by a scaling argument. This is thus the analog of the del Pino-Dolbeaut family in the trace case. 

\medskip

The paper is organized as follows. In the next section we state and prove the main results, namely generalizations of Theorem~\ref{thm-IC2} and~\ref{thm-BBL2}. In Section~\ref{sec-GN} we show how these results imply the Sobolev type inequalities: in Section~\ref{sec-rncase} we propose a new proof of the  Gagliardo-Nirenberg-Sobolev inequalities, including and extending the del Pino-Dolbeault family, whereas in Section~\ref{eq-trace-inequalities}, we follow the same procedure to reach  Gagliardo-Nirenberg-Sobolev trace inequalities. Section~\ref{sec-remarques} is devoted to classical geometric inequalities such as the Pr\'ekopa-Leindler or the classical Borell-Brascamp-Lieb inequalities, with an application to a trace logarithmic Sobolev inequality. Finally Appendix~\ref{sec-appendix} deals with a general result on the infimum convolution, which is a crucial tool for our proofs.

\medskip

Classical inequalities such as Gagliardo-Nirenberg-Sobolev are valid in $\R^n$ with some restriction on the dimension $n$. For each result,  the dimension~$n$ will be specified.

\medskip

{\bf Notation:}  When the measure is not mentionned, an integral is understood with respect to Lebesgue measure. In $\R^n$, for any $x,y\in\R^n$, $|x|$ denotes the Euclidean norm of $x$ and  $x\cdot y$ the Euclidean scalar product. As already used, $\| f\|_p$ stands for the $L^p(\R^n)$ norm.  

\bigskip

\noindent{\bf Acknowledgements.} We warmly thank S. Zugmeyer who removed an assumption on the function $\Phi$ in Theorem~\ref{thm-BBLphi}.  
 This work was partly written while the authors were visiting Institut Mittag-Leffler in Stockholm; it is a pleasure for them to thank this institution for its kind hospitality and participants for discussions on this and related works.

 This research was supported  by the French ANR-12-BS01-0019 STAB project. The third author is supported in part by JSPS KAKENHI \# 15K04949.


\section{Main results and proofs}
\label{sec-proof}

Each result presented in this paper has two formulations : the first one as a convex (or concave) inequality illustrated  by Theorem~\ref{thm-IC2} and the second one as a Borell-Brascamp-Lieb  type inequality like Theorem~\ref{thm-BBL2}. 


\subsection{Setting and additional tools}

To explain this and state our result, let us first fix the setting we are going to work with. It has two separate cases, the origin of which will be explained below. We are going to measure the gradient using a function $W$ on $\R^n$ that will belong to one of the following two categories:
\begin{enumerate}
\item Either $W:\R^n\to \R\cup\{+\infty\}$ is a {\it convex}  fonction. 
We shall let $W^\ast$ denote its Legendre transform, 
$$
W^\ast (y) = \sup_{x\in \R^n}\{ x\cdot y - W(x)\}.
$$
For almost every $x$ in the domain of $W$ the function $W$ is differentiable at $x$ and one has
\be\label{Lduality1}
W^\ast(\nabla W(x)) + W(x) = x \cdot \nabla W(x).
\ee
\item Either $W$ is a nonnegative function that is \emph{concave} on its support $\Omega_W =\{W>0\}$. More precisely, $W$ is a nonnegative  function such that the function $\tilde W$ defined on $\R^n$ by $\tilde W(x) = W(x) $ if $x\in \Omega_W$ et $-\infty$  otherwise, is concave. In particular $\Omega_W$ is a convex set. The corresponding Legendre transform is defined by 
\begin{equation}
\label{eq-inf-legendre}
W_{\ast}(y)=\inf_{x\in \Omega_W} \{x\cdot y - W(x)) \}   = \inf_{x\in \R^n} \{x\cdot y - \tilde{W}(x) \} \in \R.
\end{equation}
As above, $W$ is differentiable at almost every $x\in \Omega_W$ with 
\be \label{duality2}
W_\ast(\nabla W(x)) + W(x) = x \cdot \nabla W(x).
\ee
\end{enumerate}

We refer to the classic book  \cite{rockafellar} by R.~T.~Rockafellar for these classical definitions.\\

One rather naturally comes to such a setting if one keeps in mind the Brunn-Minkowski theory of convex measures on $\R^n$ as put forward by C. Borell. Although we will not explicitly use it, we feel it is necessary to  briefly recall it to put our results in perspective.  A nonnegative function $G$ on $\R^n$ is said to be \emph{$\kappa$-concave} with $\kappa \in \mathbb R$ if $\kappa\, G^{\kappa}$ is concave on its support. In other words, the definition splits into two categories: 
\begin{enumerate}
\item If $\kappa <0$, then $G= W^{1/\kappa}$ with $W$ convex on $\R^n$. The Brunn-Minkowski-Borell theory shows that one should consider the range $\kappa \in [-\frac1n , 0)$. Below we shall let $\kappa=-1/a$ for $a\geq n$ with the typical exemples $W(x) = 1+ |x|^q, q \geq 1$ and then $G(x)=(1+|x|^q)^{-a}$. The results above in Theorems~\ref{thm-IC2} and~\ref{thm-BBL2} correspond to the extremal case $a=n$. 
\item If $\kappa >0$, $G= W^{1/\kappa}$ with $W$ concave on its support.  Below we shall let $\kappa=1/a$ for $a>0$ with the typical examples $W(x) = (1-|x|^q)_+, q \geq 1$ and $G(x) = (1-|x|^q)_+^{a}$. 
\end{enumerate}
The limit case $\kappa = 0$ is  defined as the log-concavity of $G$. 

\bigskip

Crucial arguments in our proof are optimal transportation tools (including Brenier's map).  So let us briefly describe the mathematical setting and notation on this topic we shall use below.

\medskip

We let $\mathcal P_2(\R^n)$ be the space of probability measures $\mu$ in $\R^n$ with a finite second moment, that is $\int |x|^2d\mu(x) <+\infty$. On the optimal transportation side, Brenier's Theorem~\cite{brenier91} is the cornerstone of many proofs of functional inequalities.  It says that for any probability measure $\mu$ and $\nu$ in $\mathcal P_2(\R^n)$ with $\mu$ absolutely continuous with respect to Lebesgue measure then there exists a convex function $\phi$ (the so-called  Brenier map) on $\R^n$ such that $\nu$ is the image measure $\nabla\phi\#\mu$ of $\mu$ by $\nabla\phi$,  i.e.
 for any bounded function $H$ on $\R^n$, 
 $$
 \int Hd\nu=\int H(\nabla \phi)d\mu.
 $$ 
 From the map $\phi$ one can define a displacement interpolation from $\mu$ to $\nu$, introduced by McCann in~\cite{mccann-advances}, that is, the path $(\mu_t)_{t\in[0,1]}$ in $\mathcal P_2(\R^n)$ defined by $\mu_t=((1-t) Id +t\nabla \phi)\#\mu$, i.e.  for any bounded function $H$ 
\begin{equation}
\label{eq-geo-path}
\int Hd\mu_t=\int H((1-t)x+t\nabla \phi(x))d\mu (x). 
\end{equation}
It is now classical that  Brenier's map gives a value of the Wasserstein distance between $\mu$ and $\nu$ and $(\mu_t)_{t \in [0, 1]}$ is the geodesic in the Wasserstein space between $\mu$ an $\nu$.  These facts will not be used in our paper.

Assuming that $d\mu=fdx$ and $d\nu=gdx$ then~\cite{mccann-advances} ensures that $fdx$-almost surely, the Monge-Amp\`ere equation holds:
\begin{equation}
\label{eq-monge-ampere}
f(x)=g(\nabla \phi (x))\det(\nabla^2\phi (x)).
\end{equation}
Here $\nabla^2\phi$ is the Alexandrov Hessian of $\phi$, which is the absolutely continuous part of the distributional Hessian of the convex function $\phi$. Below we shall let $\Delta \phi$ be the trace of $\nabla^2\phi$. All these notions are explained in full details in~\cite{villani-03,villani-book1} for instance.

\medskip

Finally, our last tool will be convexity on  the determinant of matrices which we recall now (see~\cite{convex} for instance).
\begin{lemma}[Classical inequalities on the determinant]
~

\label{lem-conv}
\begin{itemize}
\item For every $k\in (0, 1/n]$, the map $H\to \det^k H$ is concave over the set of positive symmetric matrices. Concavity inequality around the identity implies 
$$
{\det}^k H \leq 1 - nk + k \tr H
$$
for all positive symmetric matrix $H$. 
\item For every $k<0$, the map $H\to \det^k H$ is convex over the set of positive symmetric matrices. Convexity inequality around the identity implies
$$
{\det}^k H \geq 1 - nk + k \tr H
$$
for all positive symmetric matrix $H$. 
\end{itemize}
\end{lemma}


\subsection{Convex and concave inequalities (Generalization of Th.~\ref{thm-IC2})}

The next two results are called {\it convex} and {\it concave} inequalities since extremal functions are convex in the first case and concave in the second one. 

\begin{theo}[Convex inequalities]
\label{thm-case1}
Let $n\geq 1$. Let $a\geq n$ (and $a>1$ if $n=1$) and let $W:\R^n \to (0,+\infty)$ such that $\int W^{1-a}<+\infty$. Then for any positive and smooth function $g$ such that $  \int W^\ast \big( \nabla g \big)g^{-a}<+\infty$, $\int g^{1-a}<+\infty$ and 
$$
\int g^{-a} = \int W^{-a}=1
$$
one has
\begin{equation}
\label{eq-case1}
(a-1) \int W^\ast \big( \nabla g \big)g^{-a} +   (a-n)\int g^{1-a}  \ge   \int W^{1-a},
\end{equation}
with equality if $g=W$ and is  convex. 
\end{theo}


\begin{eproof}
Let $\phi$ be Brenier's map such that $\nabla \phi\#g^{-a}=W^{-a}$. Then, from~\eqref{eq-monge-ampere},  almost surely,  
$$
W(\nabla \phi) = g \;  \big( \det \nabla^2\phi \big)^{1/a}. 
$$
Moreover, since $a\geq n$, from case one in Lemma~\ref{lem-conv} with $k = 1/a$ we have almost surely
$$
\big( \det \nabla^2\phi \big)^{1/a} \le 1-\frac na + \frac 1a \, \Delta \phi.
$$
Integrating with respect to the measure $g^{-a}dx$ we get 
$$
\int W(\nabla \phi) g^{-a}   \le \PAR{1-\frac na}\int g^{1-a} + \frac 1a \int  \Delta \phi \,g^{1-a},
$$
that is,
$$
a\int W(\nabla \phi) g^{-a}   \le \PAR{a-n}\int g^{1-a} + \PAR{a-1} \int  \nabla g\cdot \nabla \phi \,g^{-a}
$$
by integration by parts, justified as in~\cite[Lem.~7]{cnv04}. But
$$
\nabla g \cdot \nabla \phi \leq W(\nabla \phi) + W^*(\nabla g)
$$
almost everywhere, so collecting terms we have 
$$
\int W(\nabla \phi) g^{-a}   \le (a-1) \int  W^*(\nabla g) g^{-a}   + (a-n)\int g^{1-a}.
$$
Finally $\int W(\nabla \phi) g^{-a}=\int W^{1-a}$ since $\nabla \phi\#g^{-a}=W^{-a}$. This ends the proof of the inequality.

Now, when $g=W$ and is convex, then the relation~\eqref{Lduality1} and integration by parts ensures that inequality~\eqref{eq-case1} is an equality. 
\end{eproof}

The companion ``concave" case is as follows.

\begin{theo}[Concave inequalities]
\label{thm-case2}
Let $n\geq 1, a>0$, and $W:\R^n \to [0,+\infty)$. Then for any nonnegative  smooth function $g$ such that   
$$
\int g^{a} = \int W^{a}=1
$$
we have
\begin{equation}
\label{eq-case2}
(a+1) \int W_\ast \big( \nabla g \big)g^{a}   +  (a+n)\int g^{1+a}  \le -  \int W^{1+a}, 
\end{equation}
with equality if $g= W$ and is concave on its support (in the sense above).
\end{theo}

\begin{eproof}
The proof follows the previous one. Let $\phi$ be Brenier's map such that $\nabla \phi\#g^{a}=W^{a}$. Then, from case two in Lemma~\ref{lem-conv},
$$
W(\nabla\phi)   = g\, \PAR{\det \nabla^2\phi}^{-1/a} \ge \PAR{1+\frac na}g - \frac 1ag\, \Delta \varphi.
$$
We obtain inequality~\eqref{eq-case2} again by integrating with respect to the measure $g^adx$, integrating by parts and using the almost everywhere inequality 
$$
\nabla \phi \cdot\nabla g  \geq W(\nabla \phi) + W_\ast(\nabla g).
$$ 
When $g=W$ and is concave on its support, the inequality is an equality by~\eqref{duality2}. 
\end{eproof}


\subsection{Generalization of the Borell-Brascamp-Lieb inequality}

If Theorems~\ref{thm-case1} and~\ref{thm-case2} appear as convex or concave generalizations of Theorem~\ref{thm-IC2} (which is Theorem~\ref{thm-case1} for $a=n$), we now present two generalizations of Theorem~\ref{thm-BBL2} in the sense of Borell-Brascamp-Lieb type inequality. 

The first one concerns the convex case. 

\begin{theo}[$\Phi$-Borell-Brascamp-Lieb inequality]
\label{thm-BBLphi} 
Let $a \geq n\geq 1$ (and $a>1$ if $n=1$) and let $\Phi:\R^+\to\R^+$ be a $\mathcal C^1$-concave function.

Let also $g,W, H :\R^n \to (0,+\infty]$ be Borel functions and $t\in [0,1]$, $s=1-t$, be such that
\begin{equation}
\label{eq-contrainte}
\forall x,y\in \R^n,\qquad  H(sx+ ty) \le s g(x) + t W(y) 
\end{equation}
and $\int W^{-a} = \int g^{-a}=1$.
Then
\begin{equation}
\label{eq-BBLphi}
\int \Phi(H)H^{-a} \ge s \int \Phi(g)g^{-a} + t \int  \Phi(W)W^{-a}.
\end{equation}
\end{theo}

\begin{eproof}
The theorem can be proved in two ways, following the ideas from F. Barthe or R. J. McCann's PhDs~\cite{barthe-these,mccann-these}.

Let $\phi$ be Brenier's map such that $\nabla \phi\#g^{-a}=W^{-a}$. Then  from the Monge-Amp\`ere equation~\eqref{eq-monge-ampere},  we have that almost surely 
$$
W(\nabla \phi) = g \; \det(\nabla^2\phi)^{1/a}. 
$$

Moreover, it follows from the assumptions that $\Phi$ is nondecreasing and then  $\R^+\ni x\mapsto \Phi(x)x^{-a}$ is nonincreasing.

\medskip

\noindent
{\bf First proof:} This proof is a little bit formal since we use a change of variables formula without proof. However, it is useful to fix the ideas, and helps to follow the rigorous proof below. 

So, by change of variable and using both assumptions on $\Phi$ we have 
$$
\begin{array}{rl}
\disp\int \Phi(H)H^{-a} &\disp=\int \Phi(H(sx+t\nabla \phi(x)))H^{-a}(sx+t\nabla \phi(x))\det (s{\rm Id}+t\nabla^2\phi(x))dx\\
&\disp\geq \int {\Phi(sg+tW(\nabla \phi))} \big( sg + t W (\nabla \phi) \big)^{-a} \det (s{\rm Id}+t\nabla^2\phi).\\
&\disp\geq \int {\big[s\Phi(g)+t\Phi(W(\nabla \phi))\big]} \PAR{s+t\det(\nabla^2\phi)^{1/a}}^{-a} \det (s{\rm Id}+t\nabla^2\phi) \, g^{-a}.\\
\end{array}
$$
Since $a\geq n$, the first case in Lemma~\ref{lem-conv} with $k = 1/a$ yields 
\begin{equation}\label{eq:det}
\det (s{\rm Id}+t\nabla^2\phi)\geq \PAR{s+t\det (\nabla^2\phi)^{1/a}}^{a}.
\end{equation}
Finally $\int \Phi( W(\nabla \phi)) g^{-a}=\int \Phi(W) W^{-a}$ by image measure property since $\nabla \phi\#g^{-a}=W^{-a}$. This concludes the argument, as
$$
\disp\int \Phi(H)H^{-a} \geq \int {\big[s\Phi(g)+t\Phi(W(\nabla \phi))\big]}g^{-a}
=s\int \Phi(g)g^{-a}+t\int \Phi(W)W^{-a}.
$$

\noindent
{\bf Second proof:} We use the idea of R. J. McCann. From~\cite[Lem. D.1]{mccann-these}, let $(\rho_t)_{t\in[0,1]}$ be the density of the geodesic path between $g^{-a}$ and $W^{-a}$; as defined in~\eqref{eq-geo-path}; then almost surely
$$
\rho_t(\nabla\phi_t)\leq (sg+tW(\nabla \phi))^{-a}
$$ 
where $\phi_t(x)=s\frac{|x|^2}{2}+t\phi(x)$, $x\in\R^n$.
Multiplying the inequality by $\Phi(sg+tW(\nabla \phi))$, then a.s.
$$
\Phi(sg+tW(\nabla \phi))\rho_t(\nabla \phi_t)\leq \Phi(sg+tW(\nabla \phi))(sg+tW(\nabla \phi))^{-a}. 
$$ 
Hecne, using both assumptions on $\Phi$, we get 
$$
\big[s\Phi(g)+t\Phi(W(\nabla \phi))\big]\rho_t(\nabla \phi_t)\leq \Phi(H(sx+\nabla \phi)) \, H(sx+\nabla \phi)^{-a}
=
\Phi(H(\nabla \phi_t(x)) \, H(\nabla \phi_t(x))^{-a}. 
$$ 
Now $\nabla \phi_t\circ \nabla\phi_t^\ast={\rm Id}$ by convex analysis (a.s. in $\R^n,$ see for instance~\cite[Thm.~11 (iv)]{villani-03})   and the inequality can be written as  
$$
\big[s\Phi(g(\nabla \phi_t^\ast))+t\Phi(W(\nabla \phi(\nabla \phi_t^\ast)))\big]\rho_t\leq \Phi(H)H^{-a}. 
$$ 
Then inequality~\eqref{eq-BBLphi} follows by integration, since from~\eqref{eq-geo-path}  
\begin{multline*}
\int \big[s\Phi(g(\nabla \phi_t^\ast))+t\Phi(W(\nabla \phi(\nabla \phi_t^\ast)))\big]\rho_t=s\int \Phi(g)g^{-a}+t\int \Phi(W(\nabla \phi))g^{-a}\\
=s\int \Phi(g)g^{-a}+t\int \Phi(W)W^{-a}. 
\end{multline*} 
\end{eproof}

Theorem~\ref{thm-BBL2}  is then a particular case of Theorem~\ref{thm-BBLphi} when  $\Phi$ is the identity function and $a=n$.  Roughly  speaking, there is a hierarchy between all the family of inequalities~\eqref{eq-BBLphi} and inequality~\eqref{eq-BBL2} (when $a=n$) appears as  the strongest one.

The concave inequality in Theorem~\ref{thm-case2} also has a Borell-Brascamp-Lieb formulation. We only state it for power functions $\Phi$
since the general case is less appealing.  
\begin{theo}[A concave Borell-Brascamp-Lieb inequality]
\label{thm-BBLphi-concave}
Let $n\geq 1$ and  $a>0$.
Let also $g,W, H :\R^n \to [0,+\infty)$ be Borel functions and $t\in [0,1]$ and $s=1-t$ be such that
\begin{equation}
\label{eq-contrainte2}
\forall x,y\in \R^n,\qquad  H(sx+ ty) \geq s g(x) + t W(y) 
\end{equation}
and $\int W^{a} = \int g^{a}=1.$
Then 
\begin{equation}
\label{eq-BBL-concave}
\int H^{1+a} \ge s^{n+a+1} \int g^{1+a} + s^{n+a}t \int  W^{1+a}+(n+a)s^{n+a}t\int g^{1+a}.
\end{equation}
\end{theo}

\medskip

Inequality~\eqref{eq-BBLphi} is optimal in the sense that if $g=W$ and is convex, then one can exhibit a map $H$ which depends on $s$ such that inequality~\eqref{eq-BBLphi} is an equality. This is not the case for inequality~\eqref{eq-BBL-concave} which is less powerful than~\eqref{eq-BBLphi}.  Nevertheless the linearization of~\eqref{eq-BBL-concave} (when $t$ goes to 0) becomes optimal and gives optimal Gagliardo-Nirenberg inequalities in the concave case (cf. Section~\ref{sec-concave-GN}). 

\medskip

 \begin{eproof}
We start as in the proof of Theorem~\ref{thm-BBLphi},  sticking to the first formal argument for size limitation. Let $\phi$ be Brenier's map such that $\nabla \phi\#g^{a}=W^{a}$. Then  almost surely, 
$$
g^a = W(\nabla \phi)^{a} \; \det(\nabla^2\phi). 
$$
By assumption on $\phi$ and the concavity inequality~\eqref{eq:det}  we have 
$$
\begin{array}{rl}
\disp\int H^{1+a} &\disp=\int H^{1+a}(sx+t\nabla \phi(x))\det (s{\rm Id}+t\nabla^2\phi(x))dx\\
&\disp\geq \int \PAR{sg+tW(\nabla \phi)}^{1+a}\det (s{\rm Id}+t\nabla^2\phi)\\
&\disp\geq \int \PAR{sg+tW(\nabla \phi)}^{1+a}\PAR{s+t(\det\nabla^2\phi)^{1/n}}^n.
\end{array}
$$
 Now we keep only the order zero and one terms in the Taylor expansion in $t$ of the two terms above: 
 $$
 (sg+tW(\nabla \phi)^{1+a}
 =
 (sg)^{1+a} \Big(1 + \frac{t}{s} \frac{W(\nabla \phi)}{g} \Big)^{1+a}
 \geq
  s^{1+a}g^{1+a}+(a+1)s^{a}tg^{a}W(\nabla\phi);
 $$
 $$
 \PAR{s+t(\det\nabla^2\phi)^{1/n}}^n= s^n \PAR{1+\frac{t}{s}\PAR{\frac{g}{W(\nabla\phi)}}^{a/n}}^n\geq s^n+ns^{n-1}t\PAR{\frac{g}{W(\nabla\phi)}}^{a/n}.
 $$
 Hence
$$
\int H^{1+a}\geq s^{n+a+1} \int g^{1+a}+(1+a)s^{n+a}t\int g^a W(\nabla \phi)
+ns^{n+a}t\int g^{a}W(\nabla \phi)\PAR{\frac{g}{W(\nabla \phi)}}^{\frac{n+a}{n}}.
$$
Then in the last term we apply the inequality 
$$
nX^{\frac{n+a}{n}}\geq \PAR{n+a}X-a, 
$$
for $X \geq 0$ with $X={g}/{W(\nabla \phi)}$. We obtain the desired inequality.
 \end{eproof}
 

\subsection{Dynamical formulation of generalized Borell-Brascamp-Lieb inequality}\label{subsec-dyn}
Borell-Brascamp-Lieb inequalities admit a dynamical formulation given by the largest function $H$. Consider the following inf-convolution, defined for any functions 
$W,g:\R^n\to (0,\infty]$, $h\geq0$ and $x\in\R^n$ by
\begin{equation}
\label{eq-def-hj}
Q_h^W(g)(x)=
\left\{
\begin{array}{ll}
\disp\inf_{y\in\R^n}\BRA{g(y)+hW\PAR{\frac{x-y}{h}}}& {\rm if }\,\,h>0,\\
\disp g(x)& {\rm if }\,\,h=0
\end{array}
\right.
\end{equation}
or equivalently 
$$
Q_h^W(g)(x)=\disp\inf_{z\in\R^n}\BRA{g(x-hz)+hW(z)}.
$$
Then the constraint~\eqref{eq-contrainte} implies that  the inf-convolution
$$
H(x)=sQ_{{t}/{s}}^W(g)({x}/{s}), \qquad x\in\R^n
$$
if the best function $H$ satisfying~\eqref{eq-contrainte}.
From this observation, the $\Phi$-Borell-Brascamp-Lieb inequality~\eqref{eq-BBLphi}  admits an equivalent  dynamical formulation. 
\begin{cor}[Dynamical formulation of $\Phi$-Borell-Brascamp-Lieb] 
\label{cor-dyn-BBLphi}
Let $n\geq 1$ and   $g, W: \R^n \to (0, + \infty]$ and $\Phi$ as in Theorem~\ref{thm-BBLphi}.

If $a\geq n$ and $\int g^{-a}=\int W^{-a}=1$ then for any $h\geq 0$  the $\Phi$-Borell-Brascamp-Lieb inequality~\eqref{eq-BBLphi} is equivalent to
\begin{equation}
\label{eq-dyn-BBLphi}
(1+h)^{a-n}\int \Phi\Big(\frac{1}{1+h}Q_{h}^W(g)\Big)Q_{h}^W(g)^{-a} \ge \frac{1}{1+h}\int \Phi(g)g^{-a}+\frac{h}{1+h}\int \Phi(W)W^{-a}.
\end{equation}

In particular, when $a=n$ and $\Phi(x) = x$, the extended Borell-Brascamp-Lieb inequality~\eqref{eq-BBL2}  is equivalent to  
\begin{equation}
\label{eq-dyn-BBL2}
\forall h\geq0,\qquad\int Q_h^W(g)^{1-n} \ge \int g^{1-n}+h\int W^{1-n}.
\end{equation}
We assume that $g$ is such that all quantities are well defined in the previous two inequalities.

Let us notice that when $g=W$ and is  convex then inequalities~\eqref{eq-dyn-BBLphi} and~\eqref{eq-dyn-BBL2} are equalities.
\end{cor}

We have nothing to prove since inequality~\eqref{eq-dyn-BBLphi} is only a reformulation of~\eqref{eq-BBLphi}. We only have to check the optimal cases. When $g=W$ and is convex then, from~\eqref{eq-def-hj}, 
$$
\forall x\in\R^n,\,\,\,Q_h^W(g)(x)=(1+h) W\PAR{\frac{x}{h+1}}.
$$
It follows that, in this case,  inequalities~\eqref{eq-dyn-BBLphi} and~\eqref{eq-dyn-BBL2} are equalities. 

\medskip

Inequalities~\eqref{eq-dyn-BBLphi} and~\eqref{eq-dyn-BBL2} are equalities when $h=0$. Moreover, as  explained in Appendix~\ref{sec-appendix}, generally for $h \to 0$ we observe that 
$$
Q_h^W g = g - h W^\ast(\nabla g) + o(h),
$$
so that Theorem~\ref{thm-BBLphi} admits a linearization as a convex inequality. 
With the same conditions on the function $\Phi$ as in Theorem~\ref{thm-BBLphi},  from inequality~\eqref{eq-dyn-BBLphi} we obtain 
\begin{equation}
\label{eq-ICphi}
\int {W^\ast(\nabla g)}\PAR{a \frac{\Phi(g)}{g}-\Phi'(g)}g^{-a}+\int\PAR{\big(a-n+1\big)\Phi(g)-g\Phi'(g)}g^{-a}
\geq \int \Phi(W)W^{-a},
\end{equation}
for a class of functions $g$ and $W$ (which we do not try to carefully describe for a general $\Phi$). Of course again inequality~\eqref{eq-ICphi} is optimal: equality holds when $g=W$ and is   
convex. For the case $\Phi(x)=x$, Appendix A justifies~\eqref{eq-ICphi}, starting from~\eqref{eq-dyn-BBLphi}, for $(g,W)$ in the space $\mathcal F^a$ described in Appendix~\ref{sec-appendix-rn}, 
Definition~\ref{def-admissible+}.  This is the most important case, and then we recover Theorem~\ref{thm-case1}. When $a=n$ and again  $\Phi(x)=x$ we 
recover Theorem~\ref{thm-IC2}.

The concave Borell-Brascamp-Lieb inequality~\eqref{eq-BBL-concave} also admits a dynamical formulation with the sup-convolution instead of the inf-convolution.  Consider
$W,g:\R^n\to [0,\infty)$, $h\geq0$ and $x\in\R^n$ and let
$$
R_h^W(g)(x)=
\left\{
\begin{array}{ll}
\disp\sup_{y\in\R^n}\BRA{g(y)+hW\PAR{\frac{x-y}{h}}}& {\rm if }\,\,h>0,\\
\disp g(x)& {\rm if }\,\,h=0.
\end{array}
\right.
$$ 
Then the constraint~\eqref{eq-contrainte2} implies that the best function $H$ if given by the sup-convolution, 
$$
\forall x\in\R^n, \qquad H(x)=sR_{{t}/{s}}^W(g)({x}/{s}).
$$
From this observation, the ``concave" Borell-Brascamp-Lieb inequality~\eqref{eq-BBL-concave}  admits the equivalent following dynamical formulation: if $\int W^adx=\int g^adx=1$ then  forall $h\geq 0$,
\begin{equation}
\label{eq-dyn-BBLphi-concave}
\int R_{h}^W(g)^{1+a} \ge \int g^{1+a}+h\int W^{1+a}+\PAR{n+a}h\int g^{1+a}.
\end{equation}

Similarly to the convex cas, the properties of the semigroup $(R_{h}^W)_{h\geq0}$ ensure that the derivative of~\eqref{eq-dyn-BBLphi-concave} in $h$  (at $h=0$) implies inequality~\eqref{eq-case2}. We will not give more details on this computation.


\section{Applications to sharp Euclidean inequalities}
\label{sec-GN}

The main purpose of this section is twofolds: first we will see that the results of Section~\ref{sec-proof} imply new sharp Gagliardo-Nirenberg-Sobolev inequalities in $\R^n$. Secondly, we will give the first sharp Gagliardo-Nirenberg inequalities on the half-space $\R^n_+$, which semm completely new. 

\medskip

In all this section, $||\cdot||$ denotes an arbitrary norm in $\R^n$ and  for $y\in\R^n$ we let $||y||_*=\sup_{||x||\leq 1}x\cdot y$ its dual norm.  Recall that the Legendre transform of  $x\mapsto ||x||^q/q$ (with $q>1$) is the function $y\mapsto ||x||^p_*/p$ for $1/p+1/q=1$.


\subsection{The $\R^n$ case} 
\label{sec-rncase}

A family of sharp Gagliardo-Nirenberg-Sobolev inequalities in $\R^n$  was  first proved by M. del Pino and J. Dolbeault in~\cite{delpino-dolbeault}. The family was generalized to an arbitrary norm in~\cite{cnv04} by using the mass transportation method proposed by the second author in~\cite{cordero}.

The del Pino-Dolbeault Gagliardo-Nirenberg family of inequalities (including the Sobolev inequality) is a consequence of Theorems~\ref{thm-case1}
 and~\ref{thm-case2}.   We will prove in a rather direct and easy way that our extended Borell-Brascamp-Lieb inequality~\eqref{eq-BBL2} implies the Gagliardo-
Nirenberg-Sobolev inequalities, in the known range but also a new range of parameters. 
 As recalled in the introduction, S. Bobkov and M. Ledoux~\cite{bobkov-ledoux08} have also derived the Sobolev inequality from the Brunn-Minkowski 
 inequality, but we believe that our method is more intuitive than theirs.

\medskip

\subsubsection{From Theorem~\ref{thm-case1} to convex Gagliardo-Nirenberg-Sobolev inequalities}

Let $n\geq 1$,  $a\geq n$ ($a>1$ if $n=1$) and $q>1$.  Let $W$ be defined by $W(x)=\frac{||x||^q}{q}+C$ for $x\in\R^n$,  where the constant $C>0$ is such that $\int W^{-a}=1$. Then, for any $y\in\R^n$,  $W^\ast(y)=\frac{||y||_\ast^p}{p}-C$  where $1/p+1/q=1$. 

\begin{figure}[!h]
\begin{center}
\includegraphics[width=9cm]{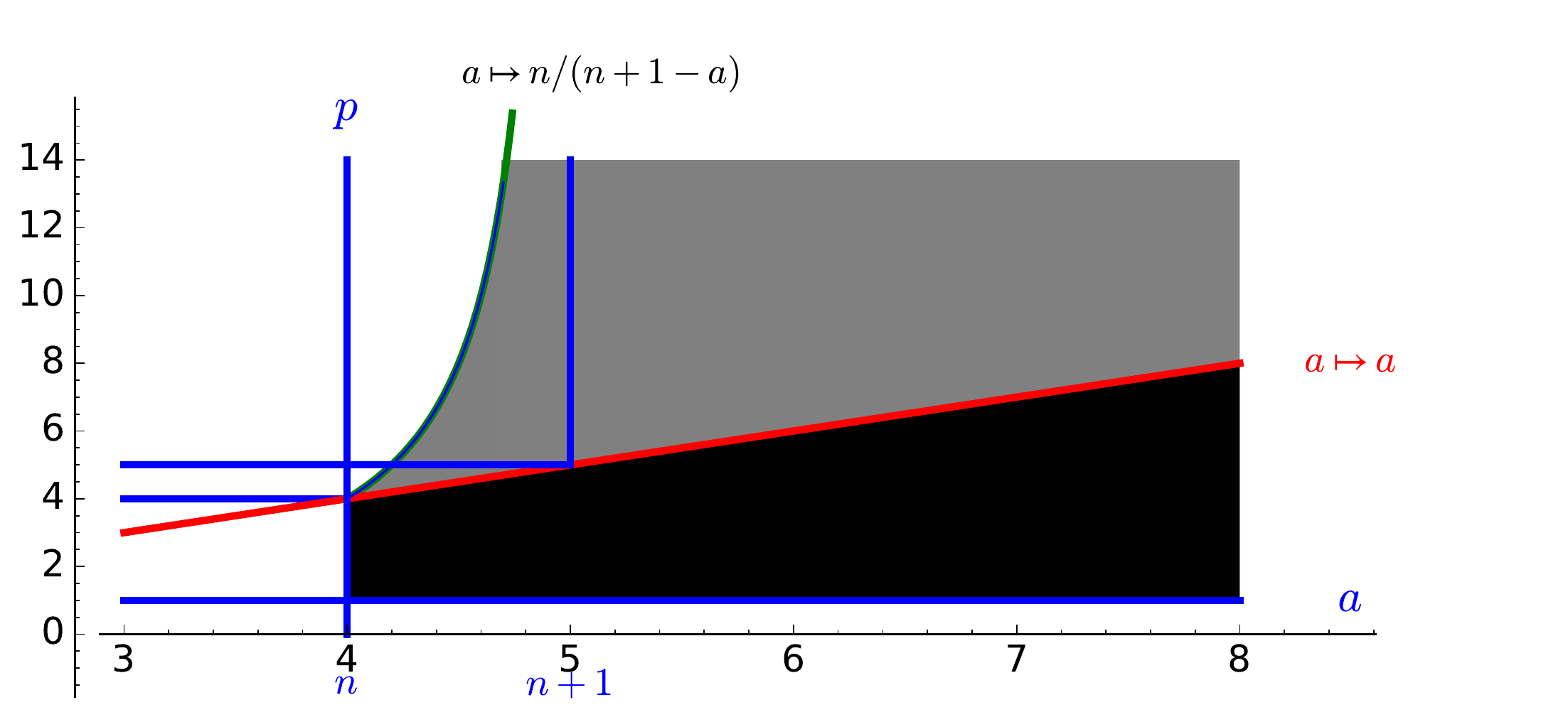}
\caption{Ranges of admissible parameters $(a,p)$ with $n=4$}
\label{fig-1}
\end{center}
\end{figure}

We would like to apply Theorem~\ref{thm-case1} with this fixed function $W$. First, let us notice that $C$ is well defined and $\int W^{1-a}$ is finite  whenever 
\begin{equation}
\label{eq-parameter}
\left\{
\begin{array}{l}
{\rm If}\,\,\, a\geq n+1\,\,\, {\rm  then}\,\,\, p>1\\ 
{\rm If}\,\,\, a\in[n,n+1)\,\,\,  {\rm then}\,\,\,   1<  p <\frac{n}{n+1-a}=\bar{p}\,\,\,  (\bar{p}=n\,\,\, {\rm  when}\,\,\,  a=n),
\end{array}
\right.
\end{equation}
These constraints are illustrated in Figure~\ref{fig-1} with the case $n=4$, Equation~\eqref{eq-parameter} is satisfied whenever the couple $(a,p)$ is in the black or the grey area.

Assuming that the parameters $a$ and $p$ are in this admissible set, then for any smooth function $g:\R^n\to \R^+$  such that $\int g^{-a}=1$, inequality~\eqref{eq-case1} in Theorem~\ref{thm-case1} becomes 
\begin{equation}
\label{eq-GN1}
D \leq  \frac{a-1}{p} \int \frac{||\nabla g||_*^p}{g^a}+(a-n)\int g^{1-a}.
\end{equation}
Here $D=(a-1)C+\int W^{1-a}$ is well defined, $W$ and $a>1$ being fixed.  This inequality is the cornerstone of this section. 

\medskip

\noindent
{\bf Sobolev inequalities: }As a warm up, let us consider $a=n$,   $n\geq 2$ and $p\in(1,n)$. Then  inequality~\eqref{eq-GN1} becomes
$$
\frac{D p}{n-1}\leq \int \frac{||\nabla g||_*^p}{g^n} 
$$
for any smooth function $g$ such that $\int g^{-n}=1$.  Letting $f=g^{\frac{p-n}{p}}$, then the inequality becomes 
$$
\frac{ D p}{n-1}\ABS{\frac{n-p}{p}}^p\leq \int {||\nabla f||_*^p}
$$
for any smooth function  $f$ such that $\int f^{\frac{np}{n-p}}=1$.  Removing the normalization we have 
$$
\frac{ D  p}{n-1}\ABS{\frac{n-p}{p}}^p\PAR{\int f^{\frac{np}{n-p}}}^{\frac{n-p}{n}}\leq \int {||\nabla f||_*^p}.
$$
The inequality is of course optimal since equality holds when $g=W$ or equivalently when $f(x)=\PAR{C+\frac{||x||^q}{q}}^{\frac{p-n}{p}}$. 
This classical result can be summarized as follows. 
\begin{theo}[Sobolev inequalities]
Let $n\geq 2$, $p\in(1,n)$ and $p^*=np/(n-p)$. The following inequality 
$$
\PAR{\int f^{p^*}}^{\frac{1}{p^*}}\leq C_{n,p}\PAR{\int {||\nabla f||_*^p}}^\frac{1}{p}.
$$
holds for any smooth function $f$ such that quantities are well defined; here $C_{n,p}$ is the optimal constant reached by the map
$\R^n\ni x\mapsto\PAR{1+{||x||^q}}^{\frac{p-n}{p}}.$
\end{theo}

\medskip

\noindent
{\bf Gagliardo-Nirenberg inequalities: }  Consider now $a > n$ (the case $a=n$ corresponds to Sobolev)  and $p\neq a$ satisfying conditions~\eqref{eq-parameter}. Letting $h=g^{\frac{p-a}{p}}$, inequality~\eqref{eq-GN1} becomes 
$$
1\leq D_2 \int {||\nabla h||_*^p}+(a-n) \int h^{p\frac{a-1}{a-p}}, 
$$
for any smooth function $h$ such that $\int h^{\frac{ap}{a-p}}=1$, where $D_2$ is an explicit positive constant.  Removing the normalization, the inequality becomes 
$$
\PAR{\int h^{\frac{ap}{a-p}}}^{\frac{a-p}{a}}
 \leq D_2 \int {||\nabla h||_*^p}+(a-n) \int h^{p\frac{a-1}{a-p}} \;  {\PAR{\int h^{\frac{ap}{a-p}}}^{\frac{1-p}{a}} } 
$$
for all smooth $h$ (such that inequalities are well defined) 

To obtain a compact form of this inequality, we replace $h(x)=f(\lambda x)$ and optimize over $\lambda >0$. We get for another explicit constant $D_3$
\begin{equation}
\label{eq-GN3}
\PAR{\int f^{\frac{ap}{a-p}}}^{\frac{a-p}{ap}\,\PAR{1-\frac{1-p}{a-p}\omega}}
\leq D_3\PAR{\int {||\nabla f||_*^p}}^{\frac{1-\omega}{p}}\PAR{\int f^{p\frac{a-1}{a-p}}}^{\frac{a-p}{p(a-1)}\frac{a-1}{a-p}\omega}
\end{equation}
where $\omega=\frac{p(a-n)}{p(a-n)+n}\in(0,1)$. There are now two cases, depending on the sign of $\PAR{1-\frac{1-p}{a-p}\,\omega}=\frac{a}{a-p}\frac{(a-n-1)p+n}{p(a-n)+n}$ 
and $\frac{a-1}{a-p}\,\omega$.  If $p<a$ then both coefficients are positive, as one can check by considering he cases $a < n+1$ and $a \geq n+1$: this leads to the first case in Theorem~\ref{thm-gn-classique} below. If $p>a$, then under the constraints~\eqref{eq-parameter} both coefficients are negative: this leads to the second case below. 

Results obtained can be summarized as follows, 
\begin{theo}[Gagliardo-Nirenberg inequalities]
\label{thm-gn-classique}
Let $n\geq 1$ and $a>n$. 

\begin{itemize}
\item For any $1<p<a$, the inequality 
\begin{equation}
\label{eq-gn-classique}
\PAR{\int f^{\frac{ap}{a-p}}}^{\frac{a-p}{ap}}\leq D_{n,p,a}^+\PAR{\int {||\nabla f||_*^p}}^{\frac{\theta}{p}}\PAR{\int f^{p\frac{a-1}{a-p}}}^{\frac{a-p}{p(a-1)}(1-\theta)}
\end{equation}
holds for any smooth function $f$ such that quantities are well defined. Here $\theta\in[0,1]$ is the unique solution of 
\begin{equation}
\label{eq-theta1}
\frac{a-p}{a} =\theta\frac{n-p}{n}+(1-\theta)\frac{a-p}{a-1}
\end{equation}
and $D_{n,p,a}^+$ is the optimal constant given by the extremal function $\R^n\ni x\mapsto \PAR{1+{||x||^q}}^{\frac{p-a}{p}}$. 

\item If $p>a$ when $a\geq n+1$, or if $p\in \big(a,\frac{n}{n+1-a} \big)$ when $a\in[n,n+1)$, then the inequality 
\begin{equation}
\label{eq-gn-classique-2}
\PAR{\int f^{p\frac{a-1}{a-p}}}^{\frac{a-p}{p(a-1)}}\leq D_{n,p,a}^-\PAR{\int {||\nabla f||_*^p}}^{\frac{\theta'}{p}}\PAR{\int f^{\frac{ap}{a-p}}}^{\frac{a-p}{ap}(1-\theta')}
\end{equation}
holds for any smooth function $f$ such that quantities are well defined. Here $\theta'\in[0,1]$ is the unique solution of 
\begin{equation}
\label{eq-theta17}
\frac{p-a}{a-1}=\theta'\frac{p-n}{n}+(1-\theta')\frac{p-a}{a} 
\end{equation}
and $D_{n,p,a}^-$ is the optimal constant given by the extremal function $\R^n\ni x\mapsto \PAR{1+{||x||^q}}^{\frac{p-a}{p}}$. In this case, the exponents in the integrals are negative. 
\end{itemize}
\end{theo}

\medskip

\begin{rem}

\label{rem-dim}
\begin{itemize}
\item Inequalities~\eqref{eq-gn-classique} is the del Pino-Dolbeault family of optimal Gagliardo-Nirenberg inequalities in $\R^n$. It correspond to parameters 
$a$ and $p$ in the black area in Figure~\ref{fig-1}.

\item   Inequalities~\eqref{eq-gn-classique-2} are Gagliardo-Nirenberg inequalities with a negative exponent, that is  
$p\frac{a-1}{a-p}<0$ and $\frac{ap}{a-p}<0$. To obtain such inequalities with the same optimal functions, the range of parameters~\eqref{eq-parameter} seems to be optimal. In this case, the couple $(a,p)$ is in the grey area in Figure~\ref{fig-1}.

Let us note that this family, with a smaller range of parameters $(a,p)$, has been obtained by V.-H. Nguyen~\cite[Th. 3.1 (ii)]{nguyen-sobolev}.
To our knowledge, the family~\eqref{eq-gn-classique-2} is new except for the part of the family proved by Nguyen. 
\item In~\cite[Th.~6.10.4]{bgl14} it has been shown how to deduce sharp Gagliardo-Nirenberg inequalities from the Sobolev inequality,  but only for the parameters $a=n+m/2$, $m\in\N$. The idea is to work in higher dimensions, for instance $\R^{n+m}$ with a function $g(x,y)=(h(x)+||y||^p)^{-(n+m-2)/2}$ and to use the scaling property of the Lebesgue measure. From inequality~\eqref{eq-case1} of Theorem~\ref{thm-case1} we can also use higher dimensions to reach all the whole family~\eqref{eq-gn-classique} of Gagliardo-Nirenberg inequalities. As in~\cite{bgl14}, we consider $g(x,y)=h(x)+||y||^r$ and $W(x,y)=||x||^p+||y||^r+C$ in $\R^{n+m}$ for a parameter $r>1$. The additional parameter $r>1$ allows us to reach all the full sharp family~\eqref{eq-gn-classique}. 
\end{itemize}
\end{rem}

\subsubsection{From Theorem~\ref{thm-case2} to concave Gagliardo-Nirenberg inequalities}
\label{sec-concave-GN}
Let $n\geq 1$.   Let $a>0$ and $q>1$, and define 
$$
\forall x\in\R^n,\,\,\, W(x)=\frac{C}{q}(1-||x||^q)_+,
$$
where $C$ is such that $\int W^a=1$. From the definition~\eqref{eq-inf-legendre}, we have 
$$
\forall y\in\R^n,\,\,\, W_*(y)=
\left\{
\begin{array}{ll}
-\frac{C^{1-p}}{p}||y||^p_*-\frac{C}{q}\,\,\, & {\rm if}\,\,\,||y||_*\leq C\\
-||y||_*\,\,\, &{\rm if}\,\,\,||y||_*\geq C
\end{array}
\right. 
$$
where $1/p+1/q=1$. 
In particular from the Young inequality 
\begin{equation}
\label{eq-maj17}
\forall y\in\R^n,\,\,\, W_*(y)\geq-\frac{C^{1-p}}{p}||y||^p_*-\frac{C}{q}.
\end{equation}
We can now apply inequality~\eqref{eq-case2} with this function $W$: for any smooth and nonnegative function $g$ such that $\int g^a=1$, 
$$
(a+n)\int g^{1+a}\leq (a+1)\frac{C^{1-p}}{p}\int {||\nabla g||_*^p}g^a+\frac{C}{q}(a+1)-\int W^{a+1}.
$$
Let us notice that $\frac{C}{q}(a+1)-\int W^{a+1}dx>0$.
Let now $h=g^{\frac{a+p}{p}}$ then, for any nonnegative function $h$ such that $\int h^{\frac{ap}{a+p}}=1$,
$$
\int h^{p\frac{1+a}{a+p}}\leq D_1\int {||\nabla h||_*^p}+D_2,
$$
where $D_1$ and $D_2$ are explicit constants. Removing the normalization, one has, for any smooth and positive function $h$,
$$
\int h^{p\frac{1+a}{a+p}}\leq D_1\int {||\nabla h||_*^p} \; \PAR{\int h^{\frac{ap}{a+p}}}^{\frac{1-p}{a}}+D_2\PAR{\int h^{\frac{ap}{a+p}}}^{\frac{1+a}{a}}.
$$
It is enough to optimize by scaling to get an inequality with a compact form. We have obtained the result proved in~\cite{delpino-dolbeault}.
\begin{theo}[Concave Gagliardo-Nirenberg inequalities]
\label{thm-gn-concave}
Let $n\geq1$.
For any $p>1$ and $a>0$ the inequality
\begin{equation}
\label{eq-gn-concave}
\PAR{\int f^{p\frac{a+1}{a+p}}}^{\frac{a+p}{p(a+1)}}\leq D_{n,p,a}\PAR{\int {||\nabla f||_*^p}}^{\frac{\theta}{p}}\PAR{\int f^{\frac{ap}{a+p}}}^{\frac{a+p}{ap}(1-\theta)}
\end{equation}
holds for any smooth and nonnegative function  $f$. Here $\theta\in[0,1]$ is the unique solution of
\begin{equation}
\label{eq-theta3}
{\frac{a+p}{a+1}} =\theta\frac{n-p}{n}+(1-\theta)\frac{a+p}{a}
\end{equation}
and  $D_{n,p,a}$ is the optimal constant given by the extremal function $\R^n\ni x\mapsto (1-||x||^q)_+^{\frac{a+p}{p}}$. 
\end{theo} 

The  obtained inequality is optimal since~\eqref{eq-case2} is an equality when $g=W$. Moreover, when $g=W$, then almost surely $||\nabla g||_*\leq C$ which implies  that~\eqref{eq-maj17} is an equality. 
\subsection{The $\R^n_+$ case, trace inequalities}
\label{eq-trace-inequalities}
For any $n\geq2$, let 
$$
\R_+^n=\{z=(u,x),\,\,u\geq0,x\in\R^{n-1}\}.
$$
Then $\partial \R^n_+=\{(0,x),\,\,x\in\R^{n-1}\}=\R^{n-1}$. For $e=(1,0)\in\R \times \R^{n-1}$ and $h\in\R$ we let 
$$
\R^{n}_{+he}=\R^{n}_++{he}=\{(u,x),\,\,u\geq h,x\in\R^{n-1}\}.
$$ 


\subsubsection{Convex inequalities in $\R^n_+$}
\label{sec-convex-trace}

The Borell-Brascamp-Lieb inequality~\eqref{eq-BBLphi} with $\Phi(x)=x$ takes the following form in $\R^n_+$.

\begin{prop}[Trace Borell-Brascamp-Lieb inequality]
Let $a\geq n$,  $g:\R^n_+ \to (0,\infty)$ and $W:\R_{+e}^n \to (0,\infty)$. Assume that $\int_{\R^n_+}g^{-a} =\int_{\R^n_{+e}}W^{-a} =1$. Then, for all $h>0$,  
\begin{equation}
\label{eq-BBL2-Trace}
(1+h)^{a-n}\int_{\R^{n}_{+he}}Q_h^W(g)^{1-a}\geq\int_{\R^n_+} g^{1-a}+h\int_{\R^n_{+e}} W^{1-a},
\end{equation}
where for any $(u,x)\in\R^{n}_{+he}$,
$$
Q_h^W(g)(u,x)=\inf_{(v,y)\in\R^n_+,\,0\leq v\leq u-h}\left\{g(v,y)+hW\PAR{\frac{u-v}{h},\frac{x-y}{h}}\right\}.
$$
Moreover~\eqref{eq-BBL2-Trace} is an equality  when $g(z)=W(z+e)$ for any $z\in\R^n_+$ and is convex. 
\end{prop}

\begin{eproof}
Let $\tilde{g}:\R^n\to (0,\infty]$ and $\tilde{W}:\R^n\to (0,\infty]$ be defined by
\begin{equation}
\label{eq-def-tilde}
\tilde{g}(x)=
\left\{
\begin{array}{l}
g(x)\,\,{\rm if}\,\,x\in\R^n_+\\
+\infty\,\,{\rm if}\,\,x\notin\R^n_+
\end{array}
\right.
\,\,\,{\rm and}\,\,\, \; 
\tilde{W}(x)=
\left\{
\begin{array}{l}
W(x)\,\,{\rm if}\,\,x\in\R^n_{+e}\\
+\infty\,\,{\rm if}\,\,x\notin\R^n_{+e}.
\end{array}
\right.
\end{equation}
Then  $\int_{\R^n}\tilde{g}^{-a}= \int_{\R^n}\tilde{W}^{-a}=1$. Hence, we can apply the dynamical formulation~\eqref{eq-dyn-BBLphi} of Theorem~\ref{thm-BBLphi} with $\Phi(x)=x$ and the functions $\tilde{g}$, $\tilde{W}$. For any $h\geq0$ we obtain 
\begin{equation*}
(1+h)^{a-n}\int_{\R^n} Q_{h}^{\tilde{W}}(\tilde{g})^{1-a} \ge \int_{\R^n_+} g^{1-a}+h\int_{\R^n_{+e}} W^{1-a},
\end{equation*}
where 
$$
\forall(u,x)\in\R^n,\quad Q_h^{\tilde{W}}(\tilde{g})(u,x)=\inf_{(v,y)\in\R^n}\left\{\tilde{g}(v,y)+h\tilde{W}\PAR{\frac{u-v}{h},\frac{x-y}{h}}\right\}.
$$
From the definition of $\tilde{g}$ and $\tilde{W}$, the infimum can be restricted to $0 \leq v \leq u-h$, so that $Q_h^{\tilde{W}}(\tilde{g})(u,x)$ is equal to $+\infty$ when $u< h,$ and to $Q_h^{{W}}({g})(x)$ otherwise.   It implies 
$$
\int_{\R^n} Q_{h}^{\tilde{W}}(\tilde{g})^{1-a}=\int_{\R^n_{+he}} Q_{h}^{\tilde{W}}(\tilde{g})^{1-a}=\int_{\R^n_{+he}} Q_{h}^{W}(g)^{1-a} ,
$$ 
which gives inequality~\eqref{eq-BBL2-Trace}. 

When $g(z)=W(z+e)$ and $W$ is convex, then by convexity 
$$
Q_h^W(g)(u,x)=(h+1)W\PAR{\frac{u+1}{h+1},\frac{x}{h+1}}
$$
for any $(u,x)\in\R^n_{+he}$. Then inequality~\eqref{eq-BBL2-Trace} is an equality.
\end{eproof}

As observed in Section~\ref{subsec-dyn}, a Borell-Brascamp-Lieb type  inequality on $\mathbb R^n$ implies a convex inequality. It is also the case on $\R^n_+$, since 
$$
\int_{\R^{n}_{+he}}Q_h^W(g)^{1-a}=\int_h^\infty\int_{\R^{n-1}}Q_h^W(g)^{1-a}(u,x)dudx
$$
and we can compute the derivative of~\eqref{eq-BBL2-Trace} on $h=0$. 

Assume that $(g,W)$ is in $\mathcal F^a_+$ as in Definition~\ref{def-admissible+}. Then by Theorem~\ref{thm-18+} in the appendix,
\begin{equation*}
\frac{d}{dh}\Big|_{h=0}\int_h^\infty\int_{\R^{n-1}}Q_h^W(g)^{1-a}(u,x)dudx
=-\int_{\partial\R^{n}_+}g^{1-a}dx+(a-1)\int_{\R^n_+}\frac{W^*(\nabla g)}{g^{a}}dz,
\end{equation*}
where we recall the definition of the Legendre transform: for any $y\in\R^n$, 
\begin{equation}
\label{def-again-legendre}
W^*(y)=\inf_{x\in\R^n_{+e}}\{x\cdot y-W(x)\}.
\end{equation}
So we have obtained:
\begin{cor}[Trace convex inequality]
\label{cor-convex-trace}
Let $a\geq n$. Let $W:\R_{+e}^n\to (0,\infty)$ be a convex function such that $\int_{\R^n_{+e}}W^{-a} =1$. Let $g:\R^n_+\to (0,\infty)$ be a smooth function satisfying  
$\int_{\R^n_+}g^{-a}=~1$.  Assume that  the couple $(g,W)$ belongs to $\mathcal F^a_+$  (see  Definition~\ref{def-admissible+} of 
Section~\ref{sec-appendix-rn+}). Then 
\begin{equation}
\label{eq-IC2-Trace}
(a-1)\int_{\R^n_+}\frac{W^*(\nabla g)}{g^{a}}+(a-n)\int_{\R^n_+}g^{1-a}\geq{\int_{\R^n_{+e}} W^{1-a}dz+\int_{\partial\R^n_+} g^{1-a}}. 
\end{equation}
Moreover~\eqref{eq-IC2-Trace} is an equality when $g(z)=W(z+e)$ for $z\in\R^n_+$ and is  convex. 
\end{cor}
\begin{rem} 
Inequality~\eqref{eq-IC2-Trace} can also be proved directly by mass transportation and integration by parts. 
\end{rem}

We follow Section~\ref{sec-rncase} to get trace Gagliardo-Nirenberg-Sobolev inequalities from Corollary~\ref{cor-convex-trace}.  
To use inequality~\eqref{eq-IC2-Trace} we need to assume that the couple $(g,W)$ is in $\mathcal F^a_+$ as in Definition~\ref{def-admissible+}. 

First we need to extend inequality~\eqref{eq-IC2-Trace} to a reasonable couple of functions $(g,W)$.

\medskip

Let $a\geq n$ and $q>1$.  Let  $W(z)=C\frac{||z||^q}{q}$ for any $z\in\R^n_{+e}$, where the constant $C>0$ is such that $\int_{\R^n_{+e}} W^{-a}=1$.  Condition {\it (C2)} of Definition~\ref{def-admissible+} is not necessarily satisfied if for instance $a\in[n,n+1)$ and $\frac{n}{a-1}\geq q>1$.  To remove this restriction we need to approximate the function $W$: for any $\varepsilon>0$ we let $\gamma>\max\{\frac{n}{a-1},1\}$ and 
$$
W_\varepsilon(x)=C_\varepsilon\frac{||z||^q}{q}+\varepsilon |x|^\gamma,
$$ 
where $C_\varepsilon$ is such that $\int_{\R^n_{+e}} W^{-a}_\varepsilon=1$.  

Then the function $W_\varepsilon$ satisfies {\it (C2)} and $\gamma$ satisfies {\it (C1)}. 
Then inequality~\eqref{eq-IC2-Trace} is valid with the function $W_\varepsilon$ for any function $g$ satisfying $\int_{\R^n_{+}} g^{-a}=1,$ {\it (C3)} and {\it (C4)}. Moreover $C_\varepsilon\rightarrow C$ and $W_\varepsilon^*(y)\rightarrow W^*(y)$, $y\in\R^n$, when $\varepsilon\rightarrow 0$.  
Then for any function $g$ satisfying $\int_{\R^n_{+}} g^{-a}=1,$ {\it (C3)} and {\it (C4)}, inequality~\eqref{eq-IC2-Trace} is valid with the function $W$. 

Now, for any $y\in\R^n$,  
\begin{equation}
\label{eq-maj1}
W^*(y)=\sup_{x\in\R^n_{+e}}\left\{x\cdot y-C\frac{||x||^q}{q}\right\}\leq\sup_{x\in\R^n}\left\{x\cdot y-C\frac{||x||^q}{q}\right\}=C^{1-p}\frac{||y||_*^p}{p}. 
\end{equation}
where $1/p+1/q=1$.  From this observation,  Corollary~\ref{cor-convex-trace} implies
\begin{equation}
\label{eq-GN4-trace}
C^{1-p}\frac{a-1}{p}\int_{\R^n_+}\frac{||\nabla g||_*^p}{g^{a}}+(a-n)\int_{\R^n_+}g^{1-a}\geq{\int_{\R^n_{+e}} W^{1-a}+\int_{\partial\R^n_+} g^{1-a}}
\end{equation}
for any function $g$ satisfying $\int_{\R^n_{+}} g^{-a}=1,$ {\it (C3)} and {\it (C4)}.

It has to be mentioned that inequality~\eqref{eq-GN4-trace} is still optimal, despite inequality~\eqref{eq-maj1}. Indeed, when $g(x)=W(x+e)$
($x\in\R^n_+$),  then the minimum in~\eqref{def-again-legendre} at the point  $\nabla g(x)$ is reached in $\R^n_{+e}$ and then~\eqref{eq-maj1} is an equality. 

Inequality~\eqref{eq-GN4-trace} is again the cornerstone of this section. 

\medskip

\noindent
{\bf Sobolev trace inequalities:} Again, as a warm up, let us assume that $a=n$. Then  the inequality~\eqref{eq-GN4-trace} becomes
$$
\int_{\partial\R^n_+} g^{1-n}\leq C^{1-p} \; \frac{n-1}{p}\int_{\R^n_+}\frac{||\nabla g||_*^p}{g^{n}}-\int_{\R^n_{+e}} W^{1-n}
$$
for any smooth function $g$ satisfying $\int g^{-n}=1$ and  {\it (C3)} and {\it (C4)}.  Assume now that $p\in(1,n)$ and let $f=g^{\frac{p-n}{p}}$. Then this inequality becomes 
\begin{equation}
\label{eq-ine-trace-encore}
\int_{\partial\R^n_+} f^{\frac{p(n-1)}{n-p}}\leq C^{1-p} \; \frac{n-1}{p}\PAR{\frac{p}{n-p}}^p\int_{\R^n_+}||\nabla f||_*^p-\int_{\R^n_{+e}} W^{1-n}
\end{equation}
under the condition $\int_{\R^n_+} f^{\frac{pn}{n-p}}=1$.  

We now need to extend the previous inequality to all smooth and compactly supported function $f$ in $\R^n_+$ (it does not mean that $f$ vanishes in $\partial\R^n_+$). For this, consider a smooth and compactly supported function $f$ in $\R^n_+$ and  let 
$$
f_\varepsilon (x)= \varepsilon |x+e|^{-\gamma \frac{n-p}{p}}+c_\varepsilon f (x),
$$
where $\gamma$ satisfies {\it (C1)} and $c_\varepsilon$ is such that $\int_{\R^n_+} f_\varepsilon^{\frac{pn}{n-p}}=1$. 
Then  
$g_\varepsilon=f_\varepsilon^{-\frac{p}{n-p}}$ satisfies {\it (C3)} and {\it (C4)}.  Moreover $c_\varepsilon\rightarrow 1$ when $\varepsilon$ goes to 0 and then
inequality~\eqref{eq-ine-trace-encore} is then valid for the function~$f$.

Removing the normalization we have for any smooth function  $f$, 
$$
\int_{\partial\R^n_+} f^{\tilde{p}}dx\leq A\int_{\R^n_+}||\nabla f||_*^pdz\,\beta^{\tilde{p}-p}-B\beta^{\tilde{p}},
$$
where
$$
\tilde{p}=\frac{p(n-1)}{n-p}, \, \, \, \beta=\PAR{\int_{\R^n_+} f^{\frac{pn}{n-p}}dz}^{\frac{n-p}{np}},\,\,\,A=C^{1-p}\frac{n-1}{p}\PAR{\frac{p}{n-p}}^p\,\,\,{\rm and}\,\,\,B=\int_{\R^n_{+e}} W^{1-n}dz.
$$
Equivalently, with  $u=\frac{\tilde{p}}{p}=\frac{n-1}{n-p}$ and $v=\frac{\tilde{p}}{\tilde{p}-p}$ (which satisfy $u,v>1$ and $1/u+1/v=1$) , 
$$
\int_{\partial\R^n_+} f^{\tilde{p}}\leq Bv\left[\frac{A}{Bv}\int_{\R^n_+}||\nabla f||_*^p\,\beta^{\tilde{p}-p}-\frac{1}{v}\beta^{\tilde{p}}\right].
$$
Now the Young inequality $xy\leq x^u/u+y^v/v$ with
$$
x=\frac{A}{Bv}\int_{\R^n_+}||\nabla f||_*^p\,\,\, \; {\rm and }\; \,\,\,y=\beta^{\tilde{p}-p}
$$
yields 
$$
\PAR{\int_{\partial\R^n_+} f^{\tilde{p}}}^{1/\tilde{p}}\leq \frac{A^{1/p}}{\PAR{Bv}^{\frac{p-1}{p(n-1)}}}\PAR{\frac{n-p}{n-1}}^{\frac{n-p}{p(n-1)}}\PAR{\int_{\R^n_+}||\nabla f||_*^p}^{1/p}.
$$


The proof of optimality it is a little bit  technical and will be given below in the more general case of Theorem~\ref{theo-GNtrace}. It is also given in~\cite{nazaret}. Equality holds when $g(z)=W(z+e)$ or equivalently when  $f(z)=\PAR{C\frac{||z+e||^q}{q}}^{-\frac{n-p}{p}}$ for $z\in\R^n_+$.

We have thus obtained the following result by B.~Nazaret~\cite{nazaret}, who promoted the idea of adding a vector $e$ to the map $W$.

\begin{theo}[Trace Sobolev inequalities from~\cite{nazaret}]\label{theo-sobtrace}
For any  $1<p<n$ and for $\tilde{p}=p(n-1)/(n-p)$ the Sobolev inequality
$$
\PAR{\int_{\partial\R^n_+} f^{\tilde{p}}dx}^{1/\tilde{p}}\leq D_{n,p}\PAR{\int_{\R^n_+}||\nabla f||_*^pdz}^{1/p}
$$
holds for any smooth function $f$ on $\R^n_+$ such that quantities are well defined. Here
$$
D_{n,p}=\frac{\PAR{\int_{\partial\R^n_+} h_p^{\tilde{p}}dx}^{1/\tilde{p}}}{\PAR{\int_{\R^n_+}||\nabla h_p||_*^pdz}^{1/p}}
$$ 
is the optimal constant, with 
$$
h_p(z)=||z+e||^{-\frac{n-p}{p-1}}, \quad z\in\R^n_+.
$$
\end{theo}

\medskip

\noindent
{\bf Gagliardo-Nirenberg trace inequalities:} Assume now that $a\geq n>p>1$ and let $h=g^{\frac{p-a}{p}}$. Then the inequality~\eqref{eq-GN4-trace} can be written as 
$$
\int_{\partial\R^n_+} h^{p\frac{a-1}{a-p}}dx\leq C^{1-p}\frac{a-1}{p}\PAR{\frac{p}{a-p}}^p\int_{\R^n_+}||\nabla h||_*^pdz+(a-n)\int_{\R^n_+}h^{p\frac{a-1}{a-p}} dz -\int_{\R^n_{+e}} W^{1-a}dz
$$
for any smooth and compactly supported functions $h$ in $\R^n_+$ such that $\int_{\R^n_+} h^{\frac{ap}{a-p}}=1$. In this case we use the same trick as for the Sobolev inequality to remove the conditions {\it (C3)} and {\it (C4)} of Definition~\ref{def-admissible+}.

Removing the normalization, then for all smooth function $h$, 
\begin{equation}
\label{eq-GN17}
\int_{\partial\R^n_+} h^{p\frac{a-1}{a-p}} \leq
C^{1-p}\frac{a-1}{p}\PAR{\frac{p}{a-p}}^p\int_{\R^n_+}||\nabla h||_*^p \,\beta^{p\frac{p-1}{a-p}} -\int_{\R^n_{+e}} W^{1-a} \,\beta^{p\frac{a-1}{a-p}}+(a-n)\int_{\R^n_+}h^{p\frac{a-1}{a-p}}
\end{equation}
with now
$$
\beta=\PAR{\int_{\R^n_+} h^{\frac{pa}{a-p}}dz}^{\frac{a-p}{ap}}.
$$
Let $u=\frac{a-1}{a-p}$ and $v=\frac{a-1}{p-1}$, which satisfy $u,v>1$ and $1/u+1/v=1$. As for the Sobolev inequality we rewrite the right-hand side of~\eqref{eq-GN17} as 
\begin{multline*}
C^{1-p}\frac{a-1}{p}\PAR{\frac{p}{a-p}}^p\int_{\R^n_+}||\nabla h||_*^p\,\beta^{p\frac{p-1}{a-p}} -\int_{\R^n_{+e}} W^{1-a}\,\beta^{p\frac{a-1}{a-p}}
\\
=
Bv\left[\frac{A}{Bv}\int_{\R^n_+}||\nabla h||_*^p\,\beta^{p\frac{p-1}{a-p}}-\frac{1}{v}\beta^{p\frac{a-1}{a-p}}\right],
\end{multline*}
with 
$$
A=C^{1-p}\frac{a-1}{p}\PAR{\frac{p}{a-p}}^p\,\,\,{\rm and}\,\,\,B=\int_{\R^n_{+e}} W^{1-a}.
$$
From the Young inequality applied to the parameters $u,v$ and
\begin{equation}
\label{xy}
x=\frac{A}{Bv}\int_{\R^n_+}||\nabla h||_*^p\,\,\, \; {\rm and }\; \,\,\,y=\beta^{p \frac{p-1}{a-p}}
\end{equation}
 we get 
\begin{equation}
\label{eq-young-trace}
C^{1-p}\frac{a-1}{p}\PAR{\frac{p}{a-p}}^p\int_{\R^n_+}||\nabla h||_*^p\,\beta^{p\frac{p-1}{a-p}} -\int_{\R^n_{+e}} W^{1-a}\,\beta^{p\frac{a-1}{a-p}}\leq
\frac{A^{\frac{a-1}{a-p}}}{\PAR{Bv}^{\frac{p-1}{a-p}}}\frac{a-p}{a-1}\PAR{\int_{\R^n_+}||\nabla h||_*^p}^{\frac{a-1}{a-p}}
\end{equation}
and then 
\begin{equation}
\label{eq-gn-lambda-avant}
\int_{\partial\R^n_+} h^{p\frac{a-1}{a-p}}dx\leq\frac{A^{\frac{a-1}{a-p}}}{\PAR{Bv}^{\frac{p-1}{a-p}}}\frac{a-p}{a-1}\PAR{\int_{\R^n_+}||\nabla h||_*^pdz}^{\frac{a-1}{a-p}}+(a-n)\int_{\R^n_+}h^{p\frac{a-1}{a-p}} dz
\end{equation}
from~\eqref{eq-GN17}. For any $\lambda>0$, we replace $h(z)=f(\lambda z)$ for $z\in\R^n_+$. We obtain 
\begin{equation}
\label{eq-gn-lambda}
\int_{\partial\R^n_+} f^{p\frac{a-1}{a-p}}dx\leq\lambda^{\frac{(a-n)(p-1)}{a-p}}\frac{A^{\frac{a-1}{a-p}}}{\PAR{Bv}^{\frac{p-1}{a-p}}}\frac{a-p}{a-1}\PAR{\int_{\R^n_+}||\nabla f||_*^pdz}^{\frac{a-1}{a-p}}+\lambda^{-1}(a-n)\int_{\R^n_+}f^{p\frac{a-1}{a-p}} dz.
\end{equation}
Taking the infimum over $\lambda>0$ gives 
$$
\PAR{\int_{\partial\R^n_+} f^{p\frac{a-1}{a-p}}dx}^{\frac{a-p}{p(a-1)}}\leq D\PAR{\int_{\R^n_+}||\nabla f||_*^pdz}^{\frac{\theta}{p}}\PAR{\int_{\R^n_+}f^{p\frac{a-1}{a-p}} dz}^{(1-\theta)\frac{a-p}{p(a-1)}}.
$$
for an explicit constant $D$ and $\theta\in[0,1]$ being the unique parameter satisfying 
\begin{equation}
\label{eq-theta2}
\frac{n-1}{n}\frac{a-p}{a-1}=\theta\frac{n-p}{n}+(1-\theta)\frac{a-p}{a-1}.
\end{equation}

We have obtained: 

\begin{theo}[Gagliardo-Nirenberg trace inequalities]\label{theo-GNtrace}
For any  $a\geq n>p>1$, the Gagliardo-Nirenberg inequality
\begin{equation}
\label{eq-gn-trace}
\PAR{\int_{\partial\R^n_+} f^{p\frac{a-1}{a-p}}dx}^{\frac{a-p}{p(a-1)}}\leq D_{n,p,a}\PAR{\int_{\R^n_+}||\nabla f||_*^pdz}^{\frac{\theta}{p}}\PAR{\int_{\R^n_+}f^{p\frac{a-1}{a-p}}}^{(1-\theta)\frac{a-p}{p(a-1)}}
\end{equation}
holds for any smooth function $f$ on $\mathbb R^n_+$ such that quantities are well defined. Here $\theta$ is defined in~\eqref{eq-theta2} and $D_{n,p,a}$ is the optimal constant, reached when
$$
f(z)=h_p(z)=||z+e||^{-\frac{a-p}{p-1}}, \quad  z\in\R^n_+.
$$
\end{theo}

When $a=n$, then $\theta=1$ and we recover the trace Sobolev inequality of Theorem~\ref{theo-sobtrace}. 

\smallskip

\begin{eproof}
From the above computation we only have to prove that the inequality~\eqref{eq-gn-trace} is optimal. 

First, it follows from Corollary~\ref{cor-convex-trace} that inequality~\eqref{eq-GN17} is an equality when 
$$
\forall z\in\R^n_+,\,\,\,h(z)=h_p(z)={{||z+e||}}^{-\frac{a-p}{p-1}},
$$ 
the function $h_p$ does not need to be normalized. Moreover, if inequality~\eqref{eq-young-trace} is an equality, then it is also the case for~\eqref{eq-gn-lambda-avant} and then~\eqref{eq-gn-trace}. So, we only have to prove that~\eqref{eq-young-trace} is an equality when $h=h_p$, which sums up to the fact that the Young inequality is an equality. This is the case when $x=y^{v-1}$ in~\eqref{xy}, that is, 
$$
\frac{A}{Bv}\int_{\R^n_+}||\nabla h||_*^pdz=\PAR{\beta^{p\frac{p-1}{a-p}}}^{v-1},
$$ 
or equivalently 
$$
\frac{A}{Bv}\int_{\R^n_+}||\nabla h||_*^pdz=\PAR{\int_{\R^n_+}{{||z+e||}}^{-\frac{ap}{p-1}}}^{\frac{a-p}{a}}.
$$ 
Let now $\mathcal I_\alpha= \int_{\R_+^n}{||z+e||}^{-\alpha}dz$ for $\alpha>0$. Then
$$
C=\frac{p}{p-1}\mathcal I_{\frac{ap}{p-1}}^{1/a},
\quad
B=\mathcal I_{\frac{ap}{p-1}}^{\frac{1-a}{a}}\mathcal I_{p\frac{a-1}{p-1}}
\quad 
{\rm{and}}
\quad 
\int_{\R^n_+}||\nabla h||_*^pdz=\PAR{\frac{a-p}{p-1}}^p\mathcal I_{p\frac{a-1}{p-1}}
$$
from their respective definition. 
Then, from the definition of $A$, equality in the Young inequality indeed holds. This finally gives equality for the map $h$. It has to be mentioned that the case $a=n$ gives the optimality of the  trace Sobolev inequality of Theorem~\ref{theo-sobtrace}.
\end{eproof}
\begin{rem}$~$
\begin{itemize}
\item We conjecture that the function $h_p$ is the only optimal function up to dilatation and translation. 
\item It was observed in~\cite{delpino-dolbeault03} that the Euclidean logarithmic Sobolev inequality can be recovered from the classical Gagliardo-Nirenberg inequality~\eqref{eq-gn-classique} by letting $a$ go to $+\infty$. In this case, a key point is that $\theta$ in equation~\eqref{eq-theta1} goes to $0$ when $a\rightarrow\infty$. In the present case of $\R^n_+$, when $a\rightarrow+\infty$, $\theta$ in equation~\eqref{eq-theta2}, goes to $1/p$: hence the method fails in $\R^n_+$. The logarithmic Sobolev inequality in $\R^n_+$ will be studied in Section~\ref{sec-remarques}.  

\item As for the Gagliardo-Nirenberg in $\R^n$, the inequality~\eqref{eq-gn-trace} can be proved by using inequality~\eqref{eq-IC2-Trace} with $a=n$ in higher dimension, as proposed in Remark~\ref{rem-dim}.
\end{itemize}
\end{rem}

\section{Remarks on classical inequalities in this context}
\label{sec-remarques}

Let us investigate, from the previous point of view,  classical inequalities as the Borell-Brascamp-Lieb and Prékopa-Leindler inequalities.  

As for the modified Borell-Brascamp-lieb inequality, the classical inequality~\eqref{eq-BBL} admits a dynamical formulation. 

Let  $W,g:\R^n\to (0, + \infty]$ satisfying $\int g^{-n}=\int W^{-n}=1$. Then, in the notation of Section~\ref{subsec-dyn}, the classical Borell-Brascamp-Lieb inequality~\eqref{eq-BBL} is equivalent to 
\begin{equation}
\label{eq-bbl-dyn}
\forall h\geq0,\qquad\int Q_h^W(g)^{-n} \ge 1.
\end{equation}
In other words, letting $\Lambda(h)=\int Q_h^W(g)^{-n}$ for $h\geq0$,  
then $\Lambda(0)=1$ and $\Lambda(h) \geq 1$ for all $h$. More surprising, it appears that $\lim_{h\rightarrow\infty}\Lambda(h)=1$, since $Q_h^W(g)(x)=hQ_{1/h}^g(W)(x/h)$ for any $h>0$ and $x\in\R^n$.

As a consequence, using the same method as in Section~\ref{subsec-dyn}, the classical Borell-Brascamp-Lieb inequality leads to the convexity 
inequality~\eqref{eq-case1} with $a=n+1$.
\begin{cor}[\cite{BGG15}]
Let  $W:\R^n\to (0, +\infty)$ be convex and such that $\int W^{-n}=1$. Then for any positive and smooth function $g$ such that $\int g^{-n}=1$,  
\begin{equation}
\label{eq-IC}
\int_{\R^n}  {W^*\PAR{\frac{\nabla g}{(\int g^{-n})^{-1/n}}}}\frac{1}{{g^{n+1}}}\geq0.
\end{equation}
\end{cor}

As we can see from Section~\ref{sec-GN}, inequality~\eqref{eq-IC} implies the family of Gagliardo-Nirenberg inequalities only for the parameters 
$a\geq n+1$. In particular, it does not imply the Sobolev inequality as pointed out  by S. Bobkov and M. Ledoux in~\cite{bobkov-ledoux08}. 

\bigskip

The Pr\'ekopa-Leindler inequality is an infinite dimensional version of the Borell-Brascamp-Lieb inequality. It states that given $H,W,g:\R^n\to \R$, $t\in[0,1]$ and $s=1-t$  satisfying $\int e^{-g}=\int e^{-W}=1$ and 
$$
\forall x,y\in \R^n,\qquad   H((1-t)x+ ty) \le(1-t)g(x)+ t W(y),
$$
then
\begin{equation*}
\int_{\mathbb R^n} e^{-H} \ge 1.
\end{equation*}
The Pr\'ekopa-Leindler inequality also admits a dynamical formulation: for any $g$ such that $\int e^{-g}dx=1$, 
\begin{equation}
\label{eq-PL}
\forall h\geq0\qquad\int_{\mathbb R^n} e^{-\frac{1}{1+h}Q_h^W(g)} \ge (1+h)^n.
\end{equation}
Again, as for previous inequalities, it admits a linearization, recovering the general logarithmic Sobolev inequality proved by the third author in~\cite{gentil03, gentil08}:

\begin{cor}[Euclidean logarithmic Sobolev inequality]
For any convex function $W:\R^n\to (0,+\infty)$ and any smooth function $g$ on $\mathbb R^n$ such that $\int e^{-g}= \int e^{-W} = 1$,  
\begin{equation}
\label{eq-IC3}
\int_{\mathbb R^n} (g+W^*(\nabla g))e^{-g}\geq n.
\end{equation}
Moreover equality holds when $g=W$ and is  convex.
\end{cor}

Let us observe quality follows from~\eqref{Lduality1}  when $g=W$ and is  convex.

\medskip

Inequality~\eqref{eq-IC3} is equivalent to 
$$ 
\ent_{dx}(f^2)=\int_{\mathbb R^n} f^2\log \frac{f^2}{\int f^2} dx \leq  \int_{\mathbb R^n} W^*\PAR{-2\frac{\nabla f}{f}}f^2dx-n\int_{\mathbb R^n} f^2dx, 
$$
for any smooth function $f$ (without normalization condition). For instance, when $W(x)=||x||^q+C$ then after scaling optimization we get the $L^p$-Euclidean logarithmic Sobolev inequality
\begin{equation}
\label{eq-lp-logsob}
\ent_{dx}(f^p)\leq  \frac{n}{p}\int_{\mathbb R^n} f^pdx \; \log\PAR{{\mathcal L_p}\frac{\int ||\nabla f||_*^pdx}{\int f^pdx}}. 
\end{equation}
Here $1/p + 1/q = 1$ and $\mathcal L_p$ is the optimal constant. It is interesting to notice that this inequality has been first obtained in~\cite{delpino-dolbeault03} as a limit case of the Gagliardo-Nirenberg inequality~\eqref{eq-gn-classique} when $a$ goes to infinity, and then generalized in~\cite{gentil03}.

\medskip

What is remarkable is that the same computation may be performed in $\R^n_+$. Indeed, as in Section~\ref{eq-trace-inequalities}, let $W:\R_{+e}^n\to \R$ and $g:\R^n_+\to \R$ such that $\int_{\R^n_+e} e^{-W} = \int_{\R^n_+} e^{-g} = 1$, and 
define $\tilde{W}$ and $\tilde{g}$ as in~\eqref{eq-def-tilde}. Then
$$
\int_{\R^n} e^{-\frac{1}{1+h}Q_{h}^{\tilde{W}}(\tilde{g})}(z)dz=\int_h^\infty\int_{\R^{n-1}} e^{-\frac{1}{1+h}Q_h^W(g)(u,x)}dudx
$$
and inequality~\eqref{eq-PL} becomes 
$$
\int_h^\infty\int_{\R^{n-1}} e^{-\frac{1}{1+h}Q_h^W(g)(u,x)}dudx\geq (1+h)^n.
$$
Its linearization when $h$ tends to $0$ is then
\begin{equation}
\label{eq-IC3-trace}
\int_{\R^n_+}(g+W^*(\nabla g))e^{-g}\geq n+\int_{\partial \R^n_+}e^{-g}
\end{equation}
whenever the  function $g$ is in a appropriate set of functions. We will not give here more details.

As in Section~\ref{sec-convex-trace}, let now $q>1$, $||\cdot||$ be a norm in $\R^n$, and let $W(z)=C\frac{||z||^q}{q}$ for $z\in\R^n_{+e}$, where $C$ is such that $\int_{\R^n_{+e}} e^{-W}=1$. Then 
$$
\forall y\in\R^n_+,\,\,\,W^*(y)\leq C^{1-p}\frac{||y||_*^p}{p}
$$
with $1/p+1/q=1$. Let then $f$ be a posiitve function on $\R^n_+$ such that $\int_{\R^n_+} f^p=1$, and apply inequality~\eqref{eq-IC3-trace} to $g=-p\log f$. After removing the normalization  we obtain  
\begin{equation}
\label{eq-vraiment-la-derniere}
\ent_{dx}(f^p)=\int_{\R^n_+} f^p\log\frac{f^p}{\int_{\R^n_+} f^pdx}dx\leq \left(\frac{C}{p}\right)^{1-p}\int_{\R^n_+}||\nabla f||_*^pdx -n\int_{\R^n_+} f^pdx-\int_{\partial \R^n_+}f^pdx 
\end{equation}
where $C=\PAR{\disp\int_{\R_+^n}e^{-\frac{||x||^q}{q}}dx}^{q/n}$. 

Inequality~\eqref{eq-vraiment-la-derniere} is a form of a trace Logarithmic Sobolev inequality. It does not have a compact expression as does inequality~\eqref{eq-lp-logsob} in the case of $\R^n$, where the scaling optimization can be performed. Nevertheless, in $\R^n_+$, it improves upon the usual~\eqref{eq-lp-logsob} if we consider functions on $\R^n_+$.


\begin{appendix}
\section{Time derivative of the infimum-convolution}
\label{sec-appendix}

The time derivative of the Hopf-Lax formula~\eqref{eq-def-hj} has been treated in different contexts, namely for Lipschitz (as in~\cite{evans}) or bounded (as in \cite{villani-book1}) initial data.
In our case the function $g$ grows as  $|x|^p$ with $p>1$ at infinity and thus these classical results can not be applied. 
We will thus follow the method proposed  by S. Bobkov and M. Ledoux~\cite{bobkov-ledoux08}, extending it to more general functions $W$ and also to the half-space $\R^n_+$. 

We give all the details for the half-space $\R^n_+$ which is the more intricate.

\subsection{The $\R^n_+$ case}
\label{sec-appendix-rn+}

Let $a\geq n$ and let $g:\R^n_+\to(0, +\infty)$, $W:\R^n_{+e}\to (0, +\infty)$ such that  $\int_{\R^n_+} g^{-a}$ and $\int_{\R^n_{+e}} W^{-a}$ are finite. The functions $g$ and $W$ are assumed to be $\mathcal C^1$ in the interior of their respective domain of definition.  Moreover we assume that $W$ goes to infinity faster that linearly:
\begin{equation}
\label{eq-sur-linear}
\lim_{x\in\R^n_{+e},|x|\to\infty}\frac{W(x)}{|x|}=+\infty.
\end{equation}

Our objective is to give sufficient conditions such that the derivative at $h=0$ of the function 
\begin{equation}
\label{eq-def-function+}
\R^+\ni h\mapsto \int_h^\infty\int_{\R^{n-1}}Q_h^W(g)^{1-a}(u,x)dudx
\end{equation}
is equal to 
$$
-\int_{\partial\R^{n}_+}g^{1-a}dz+(a-1)\int_{\R^n_+}\frac{W^*(\nabla g)}{g^{a}}dz
$$
where 
\begin{equation}
\label{def-encore-legendre}
W^*(y)=\sup_{x\in\R^n_{+e}}\{x\cdot y-W(x)\}, \qquad y \in \mathbb R^n.
\end{equation}

For this, let us first recall  the definition of  $Q_h^Wg$: for $x\in\R^n_{+{he}}$, 
\begin{equation}
\label{eqn:101}
Q_h^Wg(x)=
\left\{
\begin{array}{ll}
\displaystyle{\inf_{y\in{\R}^n_+,\,{x-y}\in\R^n_{+he}}\left[g(y)+hW\left(\frac{x-y}{h}\right)\right]}& {\rm if }\,\,h>0,\\
\displaystyle{ g(x)}& {\rm if }\,\,h=0.
\end{array}
\right.
\end{equation}
or equivalently, for $h>0$ and $x\in\R^n_{+{he}}$,
$$
Q_h^Wg(x)
=
\inf_{z\in{\R}^n_{+e},\,x-hz\in\R^n_{+}}\{g(x-hz)+hW(z)\}
=
\inf_{z\in{\R}^n_{+he},\,x-z\in\R^n_{+}}\left\{g(x-z)+hW\left(\frac{z}{h}\right)\right\}
$$

First, we have 
\begin{lemma}
\label{lem-deri-point+}
In the above notation and assumptions, for  all $x\in\overset{\circ}{\R^n_+}$
\begin{equation}
\label{eq-deri-point+}
\frac{\partial}{\partial h}\Big|_{h=0}Q_h^Wg(x)=-W^*(\nabla g(x)).
\end{equation}
\end{lemma}
\begin{eproof}
We follow and adapt the proof proposed in~\cite{bobkov-ledoux08}. Let $x\in\overset{\circ}{\R^n_+}$ be fixed.

By definition of $Q_h^W g$, for any $z\in{\R^n_{+e}}$ and $h>0$ small enough so that $x-hz\in\R^n_+$, one has
$$
\frac{Q_h^Wg(x)-g(x)}{h}\leq \frac{g(x-hz)-g(x)}{h}+W(z).
$$
Since $g$ is $\mathcal C^1$, then for all $z\in{\R^n_{+e}}$
$$
\limsup_{h\rightarrow0}\frac{Q_h^Wg(x)-g(x)}{h}\leq -\nabla g(x)\cdot z+W(z). 
$$
Then, from the definition~\eqref{def-encore-legendre} of $W^*$, 
$$
\limsup_{h\rightarrow0}\frac{Q_h^Wg(x)-g(x)}{h}\leq -W^*(\nabla g(x)). 
$$
We now prove the converse inequality. Let 
$$
A_{x,h}=\{z\in\R_{+e}^n,\,\,hW(z)\leq g(x-he)+hW(e)\}.
$$
For a small enough $h>0$ such that  $x-he\in\R_+^n$ we have   $Q_h^Wg(x)\leq g(x-he)+hW(e) $, so
$$
Q_h^Wg(x)=\inf_{z\in A_{x,h},\,x-hz\in \R_+^n}\{g(x-hz)+hW(z)\}.
$$
Hence
\begin{eqnarray*}
\frac{Q_h^Wg(x)-g(x)}{h}&=&\inf_{z\in A_{x,h},\,x-hz\in \R_+^n}\left\{\frac{g(x-hz)-g(x)}{h}+W(z)\right\}
\\&=&\inf_{z\in A_{x,h},\,x-hz\in \R_+^n}\left\{-\nabla g(x)\cdot z+z\varepsilon_x(hz)+W(z)\right\},
\end{eqnarray*}
where $\varepsilon_x(hz)\rightarrow 0$ when $hz\rightarrow 0$. It implies
$$
\frac{Q_h^Wg(x)-g(x)}{h}\geq\inf_{z\in A_{x,h}}\left\{-\nabla g(x)\cdot z+z\varepsilon_x(hz)+W(z)\right\}.
$$
By the coercivity condition~\eqref{eq-sur-linear} on $W$ and since $g$ is locally bounded, the set $A_{x,h}$ is  bounded by a constant $C$, uniformly in $h \in (0,1)$.
In particular for every $\eta>0$, there exists $h_\eta>0$ such that for all $h\leq h_\eta$ and $z\in A_{x,h}$, $|\varepsilon_x(hz)|\leq \eta$. Moreover, for all $h\leq h_\eta$, 
\begin{multline*}
\frac{Q_h^Wg(x)-g(x)}{h}\geq\inf_{z\in A_{x,h}}\left\{-\nabla g(x)\cdot z+W(z)\right\}-C\eta\geq \inf_{z\in\R_{+e}^n}\left\{-\nabla g(x)\cdot z+W(z)\right\}-C\eta\\
= -W^*(\nabla g(x))-C\eta.
\end{multline*}
Let us take the limit when $h$ goes to 0, 
$$
\liminf_{h\rightarrow0}\frac{Q_h^Wg(x)-g(x)}{h}\geq -W^*(\nabla g(x))-C\eta.
$$
As $\eta$ is arbitrary, we finally get equality~\eqref{eq-deri-point+}.
\end{eproof}

Our assumptions on the couple $(g,W)$  are summarized in the following definition. 

\begin{defi}[the set $\mathcal F^a_+$of admissible couple in $\R^n_+$]
\label{def-admissible+}
Let $n\geq2$, $g:\R^n_+\to (0,\infty)$ and  $W:\R^n_{+e}\to (0,\infty)$. We say that the couple $(g,W)$ belongs to $\mathcal F^a_+$ with $a\geq n$ if the following four conditions are satisfied for some $\gamma$: 
\begin{itemize}
\item[(C1)] $\gamma>\max\{\frac{n}{a-1},1\}$. 
\item[(C2)]
There exists a constant $A>0$ such that $W(x)\geq A|x|^\gamma$ for all $x\in {\mathbb R}^n_{+e}$.
\item[(C3)]
There exists a constant $B>0$ such that 
$
|\nabla g(x)|\leq B (|x|^{\gamma -1}+1)
$
for all $x\in {\mathbb R}^n_+.$
\item[(C4)]
There exists a constant $C$ such that $C (|x|^{\gamma } +1) \leq g(x)$ for all $x\in\R^n_{+}$. 
\end{itemize}
\end{defi}
In the following, we let $C_j$ denote several constants which are independent of $h>0$ and $x\in {\mathbb R}^n_{+he}$, but may depend on $\gamma $, $A$, $B$.

\begin{lemma}
\label{lem-lemma11+}  
Assume {\it (C1)}$\sim${\it (C4)}. Then, we find a constant $h_1>0$ such that, for all $h\in (0,h_1)$ and $x\in{\mathbb R}^n_{+he}$
\begin{equation}
\label{eq-12+}
-C_1h(1+|x|^\gamma )\leq Q_h^Wg(x)-g(x)\leq  C_2 h(|x|^{\gamma -1}+1).
\end{equation}
\end{lemma} 

\begin{eproof}
{\bf 1.} Let us first consider the easier upper bound. For any $h>0$ and $x\in\R^n_{+he}$ then $x-he\in\R^n_{+}$, so that
$$
Q_h^Wg(x)-g(x)\leq g(x-he)-g(x)+hW(e).
$$
On the other hand, for any $x\in\R^n_+$ and $y\in\R^n$ such that $x+y\in\R^n_+$ we have from {\it (C3)}, 
\begin{multline}
\label{eq-rmq}
|g(x+y)-g(x)|\\
= \Big|\int_0^1\nabla g(x+\theta y)\cdot yd\theta\Big|\leq |y| \int_0^1|\nabla g(x+\theta y)|d\theta\leq 
 C_3 |y| (|x|^{\gamma -1}+|y|^{\gamma -1}+1).
\end{multline}
From this remark applied to $y=-he$ with $h \in (0,1)$, one gets for any $x\in\R^n_{+he}$
\begin{equation}
\label{eq-maj-02}
Q_h^Wg(x)-g(x)\leq  C_4 h (|x|^{\gamma -1}+1)+hW(e)\leq C_5 h (|x|^{\gamma -1}+1).
\end{equation}

{\bf 2.} For the lower bound, we first need some preparation. Thus, fix $h\in(0,1)$ and $x\in {\mathbb R^n_{+he}}$ arbitrarily. Let $\hat y\in {\mathbb R}^n_{+he}$ be a minimizer of the infimum convolution
$$
Q_h^Wg(x)=\inf_{y\in {\mathbb R}^n_{+he}}\ \left[g(x-y)+hW\left(\frac{y}{h}\right)\right]=g(x-\hat y)+hW\left(\frac{\hat y}{h}\right).
$$
Such a $\hat y$ surely exists by {\it (C2)} and {\it (C4)}. From~\eqref{eq-maj-02} and {\it (C2)} we have (recall that $h < 1$), 
\begin{equation}
\label{eq-maj-ensuite}
\frac{A}{h^{\gamma -1}}|\hat y|^\gamma \leq hW\left(\frac{\hat y}{h}\right)\leq g(x)- g(x-\hat y) +  C_5  (|x|^{\gamma -1}+1).
\end{equation}
From inequality~\eqref{eq-rmq},
\begin{equation}
\label{eqn:14+}
|g(x)- g(x-\hat y)| \leq C_6|\hat y|\big[ |x|^{\gamma -1}+|\hat y|^{\gamma -1}+1\big].
\end{equation}
From~\eqref{eqn:14+} and~\eqref{eq-maj-ensuite}, 
$$
\frac{A}{h^{\gamma -1}}|\hat y|^{\gamma }
\leq C_6|\hat y|(|x|^{\gamma -1}+|\hat y|^{\gamma -1}+1)+C_5(|x|^{\gamma -1}+1)
$$
Choose a small constant $0<h_1\leq1$ so that
\begin{equation}
1<\frac{A}{h_1^{\gamma -1}}-C_6.\label{eqn:15+}
\end{equation}
When $0<h<h_1$, we have
$$
\frac{|\hat y|^{\gamma }}{|\hat y|+1} \leq C_7 \left[1+|x|^{\gamma -1}\right]
$$
so that 
\begin{equation}
\label{eqn:16+}
|\hat y| \leq C_8 \left(1+|x|\right).
\end{equation}

{\bf 3.} Then, fix $h\in(0,h_1)$ and $x\in {\mathbb R^n_{+he}}$ arbitrarily, where $h_1$ is the constant defined in step 2. By the arguments in step 2, we see that
\begin{equation}
\label{eqn:17+}
Q_h^Wg(x)-g(x)=\disp\inf_{y\in\R^n_{+he},\,x-y\in\R^n_+,\,|y|\leq C_8(1+|x|)} \left[g(x-y)-g(x)+hW\left(\frac{y}{h}\right)\right].
\end{equation}
As in~\eqref{eq-rmq}, we have
\begin{equation}
\label{eqn:18+}
g(x)- g(x-y)\leq |y|\ \int_0^1 |\nabla g(x-\theta y)|d\theta.
\end{equation}
When $|y|\leq C_8(1+|x|)$ and $0<\theta <1$, we have $|x-\theta y|\leq (1+C_8)(1+ |x|)$, so that 
$|\nabla g(x-\theta y)| \leq C_9(1 + |x|^{\gamma -1})$ by {\it (C3)}, uniformly in $0<\theta <1$.
Thus, when $|y|\leq C_8(1+|x|)$, we have, by~\eqref{eqn:18+},
\begin{equation*}
g(x)- g(x-y)\leq  C_9(1+|x|^{\gamma -1})\,|y|.
\end{equation*}
Hence,  by~\eqref{eqn:17+} and {\it (C1)}, we obtain
$$
\begin{array}{ll}
Q_h^Wg(x)-g(x)
&\disp\geq\inf_{y\in\R^n_{+he},\,|y|\leq C_8(1+|x|)}\ \left[-C_9(1+|x|^{\gamma -1})\,|y|+hW\left(\frac{y}{h}\right)\right]\\
&\disp\geq\inf_{y\in\R^n_{+he},\,|y|\leq C_8(1+|x|)}\ \left[-C_9(1+|x|^{\gamma -1})\,|y|+\frac{A}{h^{\gamma -1}}|y|^{\gamma }\right]\\
&\disp\geq\inf_{y\in {\mathbb R}^n}\ \left[-C_9(1+|x|^{\gamma -1})\,|y|+\frac{A}{h^{\gamma -1}}|y|^{\gamma }\right]\\
&\disp= -C_{10} h (1+|x|^{\gamma -1})^{\frac{\gamma }{\gamma -1}}.
\end{array}
$$
The last equality is a direct computation. Therefore, we conclude that
$$
Q_h^Wg(x)-g(x)\geq -C_{11} h (1+|x|^{\gamma }).
$$
The proof is complete.
\end{eproof}


\begin{lemma}
\label{lem:13+}  
Assume (C1)$\sim$(C4). Then, we find constants $C_0, h_2>0$ such that for all $h \in (0, h_2)$ and $x \in {\mathbb R}^n_{+he}$
\begin{equation}
\label{eqn:19+}
\left|\frac{Q_h^Wg(x)^{1-a}-g(x)^{1-a}}{h}\right|\leq \frac{C_0}{1 + \vert x \vert^{\gamma (a-1)}}.
\end{equation}
\end{lemma}

\begin{eproof}
First, for any  $\alpha,\beta>0$ and $a>1$, then
\begin{equation}
\label{eq-facile+}
|\alpha^{1-a}-\beta^{1-a}|\leq (a-1)|\alpha-\beta|(\alpha^{-a}+\beta^{-a}).
\end{equation}
Indeed, if for instance $\beta>\alpha>0$, then for some $\theta\in (\alpha,\beta)$ we have
$$
\alpha^{1-a} - \beta^{1-a} 
=
(a-1)(\beta-\alpha) \theta^{-a} \leq (a-1) (\beta-\alpha) \alpha^{-a}.
$$
By inequality~\eqref{eq-facile+} and Lemma~\ref{lem-lemma11+}, we have 
\begin{eqnarray*}
\left| \frac{Q_h^Wg(x)^{1-a}-g(x)^{1-a}}{h}\right|
&\leq& 
(a-1)\left|\frac{Q_h^Wg(x)-g(x)}{h}\right|\big[Q_h^Wg(x)^{-a}+g(x)^{-a}\big]
\\
&\leq&
K_1(1+|x|^\gamma )\big[ Q_h^Wg(x)^{-a} +g(x)^{-a}\big]
\end{eqnarray*}
for all $h \in (0, h_1)$ and $x \in {\mathbb R}^n_{+he}$.

On the other hand, by {\it (C4)} and Lemma~\ref{lem-lemma11+}, we have for all $h \in (0, h_1)$ and $x \in {\mathbb R}^n_{+he}$
$$
Q_h^Wg(x)\geq g(x)- C_{1} h(1+|x|^\gamma )\geq (C-C_1h)(|x|^\gamma +1).
$$
Choose a small constant $h_3$ so that
\begin{equation}
\frac{C}{2} \leq C -C_1 h_3.\label{eqn:112+}
\end{equation}
and let $h_2 \min\{h_1,h_3\}.$ Then, for all 
\begin{equation}
\label{eqn:113+}
Q_h^Wg(x)\geq \frac{C}{2} (|x|^\gamma +1)
\end{equation}
whence, again using {\it (C4)},
$$
\left| \frac{Q_h^Wg(x)^{1-a}-(g(x))^{1-a}}{h}\right| \leq C_2 (1+|x|^\gamma )^{1-a}
$$
for all  $h \in (0, h_2)$ and $x \in {\mathbb R}^n_{+he}$. 
\end{eproof}

We can now state and prove the main result of this section:   
\begin{theo}
\label{thm-18+}
In the above notation, assume that the couple $(g,W)$ is in $\mathcal F^a_+$. Then 
\begin{equation}
\label{eq-der18+}
\frac{d}{dh}\Big|_{h=0}\int_h^\infty\int_{\R^{n-1}}Q_h^W(g)^{1-a}(u,x)dudx
=
-\int_{\partial\R^{n}_+}g^{1-a}dx+(a-1)\int_{\R^n_+}\frac{W^*(\nabla g)}{g^{a}}dz.
\end{equation}
\end{theo}

\begin{eproof}
One can write the $h$-derivative as follows:
\begin{multline*}
\frac{1}{h}\PAR{\int_h^\infty\int_{\R^{n-1}}Q_h^W(g)^{1-a}(u,x)dudx-\int_{\R^{n}_+}g^{1-a}(u,x)dudx}\\
=
\frac{1}{h}\PAR{\int_h^\infty\int_{\R^{n-1}}g^{1-a}(u,x)dudx-\int_{\R^{n}_+}g^{1-a}(u,x)dudx}\\
+
\frac{1}{h}\PAR{\int_h^\infty\int_{\R^{n-1}}Q_h^W(g)^{1-a}(u,x)dudx-\int_h^\infty\int_{\R^{n-1}}g^{1-a}(u,x)dudx}.
\end{multline*}
First 
$$
\frac{1}{h}\PAR{\int_h^\infty\int_{\R^{n-1}}g^{1-a}(u,x)dudx-\int_{\R^{n}_+}g^{1-a}(u,x)dudx}=-\frac{1}{h}\int_0^h\int_{\R^{n-1}}g^{1-a}(u,x)dudx,
$$
which goes to $-\int_{\R^{n-1}}g^{1-a}(0,x)dx=-\int_{\partial\R^{n}_+}g^{1-a}$ when $h$ goes to 0.
Secondly, 
\begin{multline}
\label{eq-lebesgue}
\frac{1}{h}\PAR{\int_h^\infty\int_{\R^{n-1}}Q_h^W(g)^{1-a}(u,x)dudx-\int_h^\infty\int_{\R^{n-1}}g^{1-a}(u,x)dudx}\\
=
\int_{\R^{n}_+}\SBRA{\frac{Q_h^W(g)^{1-a}(u,x)-g^{1-a}(u,x)}{h}}1_{u\geq h}dudx.
\end{multline}
By Lemma~\ref{lem-deri-point+} the function in the right-hand side of~\eqref{eq-lebesgue} converges pointwisely to $W^*(\nabla g) g^{-a}$ as $h \to 0.$ Moreover, since $\gamma (a-1)>n$, by Lemma~\ref{lem:13+} it is bounded uniformly in $h$ by an integrable function. Hence by the Lebesgue dominated convergence Theorem the left-hand-side of~\eqref{eq-lebesgue} converges (when $h\rightarrow 0$) to
$$
(a-1)\int_{\R^{n}_+} W^*(\nabla g) g^{-a}.
$$
The proof is complete. 
\end{eproof}

\subsection{The $\R^n$ case}
\label{sec-appendix-rn}
We only give the result and conditions for the $\R^n$ case.


We let  $g:\R^n\to(0,+\infty)$ be a $\mathcal C^1$ function and $W:\R^n\to(0, +\infty)$ such that $\int g^{-n}= \int W^{-n}=1$
and
$$
\lim_{|x|\to\infty}\frac{W(x)}{|x|}=+\infty.
$$

\begin{defi}[$\mathcal F^a$, the set of admissible couple in $\R^n$]
\label{def-admissible}
Let $g:\R^n \to (0, + \infty)$ and  $W:\R^n \to (0, + \infty)$. We say that the couple $(g,W)$ belongs to $\mathcal F^n$ with $a\geq n$ ($a>1$ if $n=1$) if the following four conditions are satisfied for some $\gamma$:
\begin{itemize}
\item[(C1)] $\gamma >\max\{\frac{n}{a-1},1\}$.
\item[(C2)] 
There exists a constant $A>0$ such that 
$
W(x)\geq A|x|^\gamma
$ 
for all $x \in {\mathbb R}^n.$
\item[(C3)]  
There exists a constant $B>0$ such that 
$
|\nabla g(x)|\leq B (|x|^{\gamma -1}+1)$
for all $x\in {\mathbb R}^n.$
\item[(C4)] 
There exist a constant $C$ such that $C (|x|^{\gamma } +1) \leq g(x)$ for all $x \in {\mathbb R}^n$.
\end{itemize}
\end{defi}

\begin{theo}
\label{thm:15}  
Assume that the couple $(g,W)$ is in $\mathcal F^a$. Then, the derivative at $h=0$ of the map 
$$
(0,\infty)\ni h\mapsto \int Q_h^W(g)^{1-a}
$$
is equal to 
$$
(1-a)\int  \frac{W^*(\nabla g)}{g^a}.
$$
\end{theo}

\end{appendix}


\begin{thebibliography}{{Aub}76}

\bibitem[{Aub}76]{aubin76}
T.~{Aubin}.
\newblock {Probl\`emes isoperimetriques et espaces de Sobolev.}
\newblock {\em {J. Differ. Geom.}}, 11:573--598, 1976.

\bibitem[Bar97]{barthe-these}
F.~Barthe.
\newblock {\em Inégalités fonctionnelles et géométriques obtenues par
  transport des mesures}.
\newblock PhD Thesis, 1997.

\bibitem[BGG15]{BGG15}
F.~Bolley, I.~Gentil, and A.~Guillin.
\newblock Dimensional improvements of the logarithmic {S}obolev, {T}alagrand
  and {B}rascamp-{L}ieb inequalities.
\newblock {\em Preprint}, 2015.

\bibitem[BGL14]{bgl14}
D.~{Bakry}, I.~{Gentil}, and M.~{Ledoux}.
\newblock {\em {Analysis and geometry of Markov diffusion operators.}}
\newblock Springer, Cham, 2014.

\bibitem[BL76]{brascamp-lieb76}
H.~J. {Brascamp} and E.~H. {Lieb}.
\newblock {On extensions of the Brunn-Minkowski and Prekopa-Leindler theorems,
  including inequalities for log concave functions, and with an application to
  the diffusion equation.}
\newblock {\em {J. Funct. Anal.}}, 22:366--389, 1976.

\bibitem[BL00]{bobkov-ledoux00}
S.G. {Bobkov} and M.~{Ledoux}.
\newblock {From Brunn-Minkowski to Brascamp-Lieb and to logarithmic Sobolev
  inequalities.}
\newblock {\em {Geom. Funct. Anal.}}, 10(5):1028--1052, 2000.

\bibitem[BL08]{bobkov-ledoux08}
S.~G. {Bobkov} and M.~{Ledoux}.
\newblock {From Brunn-Minkowski to sharp Sobolev inequalities.}
\newblock {\em {Ann. Mat. Pura Appl. (4)}}, 187(3):369--384, 2008.

\bibitem[{Bor}75]{borell75}
C.~{Borell}.
\newblock {Convex set functions in $d$-space.}
\newblock {\em {Period. Math. Hung.}}, 6:111--136, 1975.

\bibitem[{Bre}91]{brenier91}
Y.~{Brenier}.
\newblock {Polar factorization and monotone rearrangement of vector-valued
  functions.}
\newblock {\em {Commun. Pure Appl. Math.}}, 44(4):375--417, 1991.

\bibitem[BV04]{convex}
S.~Boyd and L.~Vandenberghe.
\newblock {\em Convex optimization}.
\newblock Cambridge Univ. Press, Cambridge, 2004.

\bibitem[CE02]{cordero}
D.~Cordero-Erausquin.
\newblock Some applications of mass transport to {G}aussian-type inequalities.
\newblock {\em Arch. Rational Mech. Anal.}, 161(3):257--269, 2002.

\bibitem[CNV04]{cnv04}
D.~{Cordero-Erausquin}, B.~{Nazaret}, and C.~{Villani}.
\newblock {A mass-transportation approach to sharp Sobolev and
  Gagliardo-Nirenberg inequalities.}
\newblock {\em {Adv. Math.}}, 182(2):307--332, 2004.

\bibitem[dD02]{delpino-dolbeault}
M.~{del Pino} and J.~{Dolbeault}.
\newblock {Best constants for Gagliardo-Nirenberg inequalities and applications
  to nonlinear diffusions.}
\newblock {\em {J. Math. Pures Appl. (9)}}, 81(9):847--875, 2002.

\bibitem[dD03]{delpino-dolbeault03}
M.~{del Pino} and J.~{Dolbeault}.
\newblock {The optimal Euclidean $L^{p}$-Sobolev logarithmic inequality.}
\newblock {\em {J. Funct. Anal.}}, 197(1):151--161, 2003.

\bibitem[{Eva}98]{evans}
L.~C. {Evans}.
\newblock {\em {Partial differential equations.}}
\newblock Amer. Math. Soc., Providence, 1998.

\bibitem[Gar02]{gardner02}
R.~J. Gardner.
\newblock The {B}runn-{M}inkowski inequality.
\newblock {\em Bull. Amer. Math. Soc. (N.S.)}, 39(3):355--405, 2002.

\bibitem[{Gen}03]{gentil03}
I.~{Gentil}.
\newblock {The general optimal $L^{p}$-Euclidean logarithmic Sobolev inequality
  by Hamilton--Jacobi equations.}
\newblock {\em {J. Funct. Anal.}}, 202(2):591--599, 2003.

\bibitem[{Gen}08]{gentil08}
I.~{Gentil}.
\newblock {From the Pr\'ekopa-Leindler inequality to modified logarithmic
  Sobolev inequality.}
\newblock {\em {Ann. Fac. Sci. Toulouse, Math. (6)}}, 17(2):291--308, 2008.

\bibitem[McC94]{mccann-these}
R.~J. McCann.
\newblock {\em A Convexity Theory for Interacting Gases and Equilibrium
  Crystals}.
\newblock PhD Thesis, 1994.

\bibitem[McC97]{mccann-advances}
R.~J. McCann.
\newblock A convexity principle for interacting gases.
\newblock {\em Adv. Math.}, 128:153--179, 1997.

\bibitem[{Naz}06]{nazaret}
B.~{Nazaret}.
\newblock {Best constant in Sobolev trace inequalities on the half-space.}
\newblock {\em {Nonlinear Anal., Theory Methods Appl., Ser. A, Theory
  Methods}}, 65(10):1977--1985, 2006.

\bibitem[Ngu15]{nguyen-sobolev}
V.-H. Nguyen.
\newblock {Sharp weighted Sobolev and Gagliardo-Nirenberg inequalities on
  half-spaces via mass transport and consequences.}
\newblock {\em {Proc. Lond. Math. Soc. (3)}}, 111(1):127--148, 2015.

\bibitem[Roc70]{rockafellar}
R.T. Rockafellar.
\newblock {\em Convex analysis}, volume~28 of {\em Princeton Math. Series}.
\newblock Princeton Univ. Press, Princeton, 1970.

\bibitem[{Tal}76]{talenti76}
G.~{Talenti}.
\newblock {Best constant in Sobolev inequality.}
\newblock {\em {Ann. Mat. Pura Appl. (4)}}, 110:353--372, 1976.

\bibitem[Vil03]{villani-03}
C.~Villani.
\newblock {\em Topics in optimal transportation}, volume~58 of {\em Grad.
  Studies Math.}
\newblock Amer. Math. Soc., Providence, 2003.

\bibitem[Vil09]{villani-book1}
C.~Villani.
\newblock {\em Optimal transport, Old and new}, volume 338 of {\em Grund. Math.
  Wiss.}
\newblock Springer, Berlin, 2009.

\end{thebibliography}
\end{document}